\documentclass{article}

\usepackage{theorem}
\usepackage{amssymb}
\usepackage{amsmath}
\usepackage{amsfonts}
\usepackage{times}
\usepackage{graphicx}

\def \RR {\mathbb R}

\def \EE {\mathbb E}

\def \PP {\mathbb P}
\def \eps {\varepsilon}
\def \vphi {\varphi}
\def \bM {\mathbb M}

\def \cE {\mathcal E}
\def \cF {\mathcal F}

\def \cN {\mathcal N}

\def \cI {\mathcal I}
\def \cR {\mathcal R}

\def \cO {\mathcal O}
\def \cV {\mathcal V}
\def \cG {\mathcal G}

\newtheorem{theorem}{Theorem}[section]
\newtheorem{lemma}[theorem]{Lemma}

\newtheorem{corollary}[theorem]{Corollary}
 {\theorembodyfont{\rmfamily}}
\newtheorem{definition}[theorem]{Definition}

\def\myffrac#1#2 in #3{\raise 2.6pt\hbox{$#3 #1$}\mkern-1.5mu\raise 0.8pt\hbox{$
#3/$}\mkern-1.1mu\lower 1.5pt\hbox{$#3 #2$}}

\begin{document}

\title{On nearly radial marginals
of high-dimensional probability measures}
\author{
Bo'az Klartag\thanks{Supported in part by the Israel Science
Foundation and by a Marie Curie Reintegtation Grant from the
Commission of the European Communities.}}
\date{}

\maketitle
\begin{abstract}
 Suppose that $\mu$ is an absolutely continuous probability measure
on $\RR^n$, for large $n$. Then  $\mu$ has low-dimensional marginals
that are approximately spherically-symmetric. More precisely, if $n
\geq (C / \eps)^{C d}$, then there exist $d$-dimensional marginals
of $\mu$ that are $\eps$-far from being spherically-symmetric, in an
appropriate sense. Here $C > 0$ is a universal constant.
\end{abstract}

\section{Introduction}

The purpose of this paper is to clarify a seven line paragraph by
Gromov \cite[Section 1.2.F]{gromov}. We are interested in
projections of high-dimensional probability measures. Not all
probability measures on $\RR^n$, for large $n$, are truly
$n$-dimensional. For instance, a measure supported on an atom or two
should not be considered  high-dimensional. Roughly speaking, we
think of a probability measure on a linear space as decently
high-dimensional if any subspace of bounded dimension contains only
a small fraction of the total mass.

\begin{definition} Let $\mu$ be a Borel probability measure on $\RR^n$
and $\eps > 0$. We say that $\mu$ is ``decently high-dimensional
with parameter $\eps$'', or ``$\eps$-decent'' in short, if for any
linear subspace $E \subseteq \RR^n$,
\begin{equation}
 \mu(E) \leq \eps \dim(E). \label{eq_613}
\end{equation}
We say that $\mu$ is decent if it is $\eps$-decent for $\eps = 1/n$,
the minimal possible value of $\eps$. \label{def_1254}
\end{definition}

Clearly, all absolutely continuous probability measures on $\RR^n$
are decent, as are many discrete measures. Note that a decent
measure $\mu$ necessarily satisfies $\mu (\{0 \}) = 0$, however,
this feature should not be taken too seriously. A measure $\mu$ is
``weakly $\eps$-decent'' if (\ref{eq_613}) holds for all subspaces
$E \subseteq \RR^n$ except $E = \{ 0 \}$. For a measure $\mu$ on a
measurable space $\Omega$ and a measurable map $T: \Omega
\rightarrow \Omega^{\prime}$, we write $T_*(\mu)$ for the
push-forward of $\mu$ under $T$, i.e.,
$$ T_*(\mu)(A) = \mu(T^{-1}(A)) $$
for all measurable sets $A \subseteq \Omega^{\prime}$. When $\mu$ is
a probability measure on $\RR^n$ and $T: \RR^n \rightarrow \RR^\ell$
is a {\it linear} map with $\ell < n$, we say that $T_*(\mu)$ is a
marginal of $\mu$, or a measure projection of $\mu$. The classical
Dvoretzky theorem asserts that appropriate geometric projections of
any high-dimensional convex body are approximately Euclidean balls
(see Milman \cite{mil} and references therein). The analogous
statement for probability measures should perhaps be the following
(see Gromov \cite{gromov}): Appropriate measure projections of any
decent high-dimensional probability measure are approximately
spherically-symmetric. When can we say that a probability measure
$\mu$ on $\RR^d$ is approximately radially-symmetric?

\medskip
We need some notation. Let $\mu$ be a finite measure on
a measurable space $\Omega$.
For a subset $A \subseteq \Omega$ with $\mu(A) > 0$ we write
$\mu|_A$ for the conditioning of $\mu$ on $A$, i.e.,
$$ \mu|_A(B) = \mu(A \cap B) / \mu(A) $$
for any measurable set $B \subseteq \Omega$. Write $S^{d-1}$ for the
unit sphere centered at the origin in $\RR^d$. The uniform
probability measure on the sphere $S^{d-1}$ is denoted by
$\sigma_{d-1}$. For two probability measures $\mu$ and $\nu$ on the
sphere $S^{d-1}$ and $1 \leq p < \infty$ we write $W_p(\mu, \nu)$
for the $L^p$ Monge-Kantorovich transportation distance between
$\mu$ and $\nu$ in the sphere $S^{d-1}$ endowed with the geodesic
distance (see, e.g., \cite{villani} or Section \ref{sec2} below).
The metrics $W_p$ are all equivalent (we have $W_1 \leq W_p \leq \pi
W_1^{1/p}$) and they metrize weak convergence of probability
measures. For an interval $J \subset (0, \infty)$ we consider the
spherical shell $S(J) = \{ x \in \RR^d ; |x| \in J \}$, where $|
\cdot |$ is the standard Euclidean norm in $\RR^d$. The radial
projection in $\RR^d$ is the map $\cR(x) = x / |x|$. An interval is
either open, closed or half-open and half-closed.

\begin{definition} [Gromov \cite{gromov_private}]
Let $\mu$ be a Borel probability measure on $\RR^d$
and let $\eps > 0$. We say that $\mu$ is ``$\eps$-radial'' if
for any interval $J \subset (0, \infty)$ with $\mu(S(J)) \geq \eps$,
we have
$$ W_{1} \left( \cR_*(\mu|_{S(J)}), \sigma_{d-1} \right) \leq \eps. $$
That is, when we condition $\mu$ on any spherical shell that
contains at least an $\eps$-fraction of the mass, and then project
radially to the sphere, we obtain an approximation to the uniform
probability measure on the sphere in the transportation-metric
sense. \label{def_934}
\end{definition}

Note that this definition is scale-invariant. We think of the
dimension $n$ from Definition \ref{def_1254} as a very large number,
tending to infinity. On the other hand, we usually view the
dimension $d$ in Definition \ref{def_934} as being fixed, and
typically not very large. The case $d=1$ of Definition \ref{def_934}
corresponds to the measure being approximately even. We are not sure
whether Dirac's measure $\delta_0$ is a good example of an
$\eps$-radial measure. An $\eps$-radial measure $\mu$ is said to be
``proper'' if $\mu (\{ 0 \}) = 0$. Our main theorem reads as
follows:

\begin{theorem} There exists a universal constant $C > 0$ for which the following holds: Let $0 < \eps < 1$ and let $d, n$ be  positive integers.
Suppose that \begin{equation} n \geq (C / \eps)^{C d}.
\label{eq_1521} \end{equation} Then, for any decent
 probability measure $\mu$ on $\RR^n$, there exists
  a linear map $T: \RR^n \rightarrow \RR^d$ such that $T_*(\mu)$ is $\eps$-radial proper.

\medskip Furthermore, let $\eta > 0$ be such that $\eta^{-1} \geq (C / \eps)^{C d}$. Then, for any $\eta$-decent probability measure $\mu$ on $\RR^n$, there exists a linear map $T: \RR^n \rightarrow \RR^d$ such that $T_*(\mu)$ is $\eps$-radial proper.
 \label{main_theorem}
\end{theorem}

Gromov has a topological proof for the cases $d = 1,2$ of Theorem
\ref{main_theorem}, which does not seem to generalize to higher
dimensions \cite{gromov_private}, \cite{few}. Theorem
\ref{main_theorem} is tight, up to the value of the constant $C$, as
demonstrated by the example where $\mu$ is distributed uniformly on
$n$ linearly independent vectors: In this case $\mu$ is decent, but
for any linear map $T$ and an interval $J$, the discrete measure
 $\cR( (T_* \mu)|_{S(J)})$ is composed of
 at most $n$ atoms. It is not difficult to see that when the support of $\nu$ contains
no more than $\eps^{-(d-1)}$ points, we have the lower bound
$W_1(\nu, \sigma_{d-1}) \geq c \eps$, for a certain universal
constant $c > 0$.
 It is desirable to find the best constant in the exponent
in Theorem \ref{main_theorem}, perhaps also with respect to other
notions of $\eps$-radial measures.

\medskip The conclusion of Theorem \ref{main_theorem} also holds
when the measure $\mu$ is assumed to be only ``weakly
$\eps$-decent'', except that $T_*(\mu)$ is no longer necessarily
proper. Another possibility in this context is to allow affine maps
in Theorem \ref{main_theorem} in place of linear maps, and obtain a
measure $T_*(\mu)$ which is $\eps$-radial proper. (It is also
possible to modify Definition \ref{def_1254} slightly, and require
that (\ref{eq_613}) hold  for all affine subspaces of dimension at
least one. The effect of such a modification is minor, since an
$\eps$-decent measure will remain at most $2 \eps$-decent after such
a change).

\medskip
The conclusion of Theorem \ref{main_theorem} does not necessarily
hold for non-decent measures, even when their support spans the
entire $\RR^n$: Let $e_1,\ldots,e_n$ be linearly independent
vectors in $\RR^n$, and consider the probability measure $\mu = (1 -
2^{-n})^{-1} \sum_{i=1}^n 2^{-i} \delta_{e_i}$, where $\delta_x$ is
Dirac's unit mass at $x \in \RR^n$. Then $\mu$ is not decent, and
none of the two-dimensional marginals of $\mu$ are $\eps$-radial
proper, for $\eps = 1/10$.

\medskip As in Milman's proof of Dvoretzky's theorem (see \cite{mil}),
Theorem \ref{main_theorem} will be proved by demonstrating that a
random linear map $T$ works with positive probability, once the
measure $\mu$ is put in the right ``position''. That is, we first
push-forward $\mu$ under an appropriate invertible linear map in
$\RR^n$, which is non-random, and only then do we project the
resulting probability measure to a random $d$-dimensional subspace,
distributed uniformly in the Grassmannian. The measure $\mu$ is in
the correct ``position'' when the covariance matrix of $\cR_* \mu$
is proportional to the identity matrix. If we assume that the
covariance matrix of $\mu$ itself is proportional to the identity,
then a random projection will not work, in general, with high
probability (compare with Sudakov's theorem; see \cite{sudakov} or
the presentation in Bobkov \cite{bobkov}).

\medskip
Here is an outline of the proof of Theorem \ref{main_theorem} and
also of the structure of this manuscript: In Section \ref{sec5} we
use the non-degeneracy conditions from Definition \ref{def_1254} in
order to guarantee the existence of the initial linear
transformation that puts $\mu$ in the right ``position''. Once we
know that the covariance matrix of $\cR_* \mu$ is approximately a
scalar matrix, we prove that the measure $\mu$ may be decomposed
into many almost-orthogonal ensembles. Each such ensemble is simply
a discrete probability measure, uniform on a collection of
approximately-orthogonal vectors in $\RR^n$ that are not necessarily
of the same length. This decomposition, which essentially appeared
earlier in the work of Bourgain, Lindenstrauss and Milman
\cite{BLM}, is discussed in Section \ref{sec4}.  Section \ref{sec3}
is concerned with the analysis of a single ensemble of our
decomposition. As it turns out, a random projection works with high
probability, and transforms the discrete measure into an
almost-radial one.
 Section \ref{sec2} contains a
preliminary discussion regarding $\eps$-radial measures and the
transportation metric. The proof of Theorem \ref{main_theorem} is
completed in Section \ref{sec6}, in which we also make some related
comments and prove the following corollary to Theorem
\ref{main_theorem}.

\begin{corollary}
There exists a sequence $R_n \rightarrow \infty$ with the following
property: Let $\mu$ be a decent probability measure on $\RR^n$.
Then, there exists a non-zero linear functional
$\vphi: \RR^n \rightarrow \RR$ such that
$$ \mu \left( \{ x  ; \vphi(x) \geq t M \right \}) \geq c \exp(-C t^2)
\ \ \ \ \ \ \text{for all} \ 0 \leq t \leq R_n $$ and
$$ \mu \left( \{ x  ; \vphi(x) \leq -t M \right \}) \geq c \exp(-C t^2)
\ \ \ \ \ \ \text{for all} \ 0 \leq t \leq R_n $$ where $M > 0$ is a
median, that is,
\begin{equation}
 \mu \left( \{ x  ; |\vphi(x)| \leq M \right \})
\geq 1/2 \ \ \ \ \text{and} \ \ \ \ \mu \left( \{ x  ; |\vphi(x)|
\geq M \right \}) \geq 1/2 \label{eq_1159} \end{equation} and $c, C
> 0$ are universal constants. Moreover, one may take $R_n = c (\log
n)^{1/4}$. \label{cor_1222}
\end{corollary}

In other words, any high-dimensional probability measure has
super-gaussian marginals. Furthermore, as is evident from the proof,
{\it most} of the marginals are super-gaussian when the measure is
in the right ``position''. In the case of independent random
variables, Corollary \ref{cor_1222} essentially goes back to
Kolmogorov \cite{kol}. See also Nagaev \cite{nag}.

\medskip In Section \ref{sec7} we formulate our results in an
infinite-dimensional setting. Unless stated otherwise, throughout
the text the letters $c, C, C^{\prime}, \tilde{c}$ etc. stand for
various positive universal constants, whose value may change from
one instance to the next. We usually denote by lower-case $c,
\tilde{c}, c^{\prime}, \bar{c}$ etc. positive universal constants
that are assumed sufficiently small, and by upper-case $C,
\tilde{C}, C^{\prime}, \bar{C}$ etc. sufficiently large universal
constants. We write $x \cdot y$ for the usual scalar product of $x,y
\in \RR^n$.

\medskip
\emph{Acknowledgments.} I would like to thank Misha Gromov for
 his interest in this work and for introducing me to the problem,
to Vitali Milman for encouraging me to write this note, to Boris
Tsirelson for his explanations regarding measures on infinite-dimensional
linear spaces, to Sasha Sodin for reading a preliminary version of
this text and to Semyon Alesker, Noga Alon and Apostolos Giannopoulos for
related discussions.

\section{Transportation distance and empirical distributions}
\label{sec2}

Let $(X, \rho)$ be a metric space and let $\mu_1, \mu_2$ be Borel
probability measures on $X$. A {\it coupling} of $\mu_1$ and $\mu_2$
is a Borel probability measure $\gamma$ on $X \times X$ whose first
marginal is $\mu_1$ and whose second marginal is $\mu_2$, that is,
$(P_1)_* \gamma = \mu_1$ and $(P_2)_* \gamma = \mu_2$ where
$P_1(x,y) = x$ and $P_2(x,y) = y$. The $L^1$ Monge-Kantorovich
distance is $$ W_{1}(\mu_1,\mu_2) = \inf_{\gamma} \int_{X \times X}
\rho(x,y) d \gamma(x,y) $$ where the infimum runs over all couplings
$\gamma$ of $\mu_1$ and $\mu_2$. Then $W_1$ is a metric, and it satisfies
the convexity relation
\begin{equation}
 W_1 \left( \lambda \mu_1 + (1 - \lambda) \mu_2, \nu \right)
\leq \lambda W_1(\mu_1, \nu) + (1 - \lambda) W_1(\mu_2, \nu)
\label{eq_1025}
\end{equation}
for any $0 < \lambda < 1$ and probability measures $\mu_1, \mu_2,
\nu$ on $X$. The Kantorovich-Rubinstein duality theorem (see
\cite[Theorem 1.14]{villani}) states that
\begin{equation}
W_1(\mu, \nu) = \sup_{\vphi}  \int_{X} \vphi d \left[ \mu - \nu \right] \label{eq_1455}
\end{equation}
where the supremum runs over all $1$-Lipschitz functions $\vphi:X
\rightarrow \RR$ (i.e., $|\vphi(x) - \vphi(y)| \leq \rho(x,y)$ for
all $x,y \in X$).  We are concerned mostly with the case where the
metric space $X$ is the Euclidean sphere $S^{n-1}$ with the metric
$\rho(x,y)$ being the geodesic distance in $S^{d-1}$, i.e., $\cos
\rho(x,y) = x \cdot y$. Denote by $\bM(S^{d-1})$ the space of Borel
probability measures on $S^{d-1}$, endowed with the weak$^*$
topology and the corresponding Borel $\sigma$-algebra. Similarly,
$\bM(\RR^d)$ is the space of Borel probability measures on
$\RR^{d}$, endowed with the weak$^*$ topology (convergence of
integrals of compactly-supported continuous functions) and
$\sigma$-algebra. A measure here always means a non-negative
measure. The total variation distance between two measures $\mu$ and
$\nu$ on a measurable space $\Omega$ is
$$ d_{TV}(\mu, \nu) = \sup_{A \subseteq \Omega} \left|
\mu(A) - \nu(A) \right| $$ where the supremum runs over all
measurable sets $A \subseteq \Omega$. Clearly, for $\mu, \nu \in
\bM(S^{d-1})$,
\begin{equation}
W_1(\mu, \nu) \leq \pi d_{TV}(\mu, \nu) \leq \pi.
\label{eq_925}
\end{equation}
Additionally, $d_{TV}(S_* \mu, S_* \nu) \leq d_{TV}(\mu, \nu)$ for
any measures $\mu, \nu$ and a measurable map $S$. When $S$ is a
$\lambda$-Lipschitz map between metric spaces, we obtain the
inequality $W_1(S_* \mu, S_* \nu) \leq \lambda W_1(\mu, \nu)$. The
following lemma is an obvious consequence of (\ref{eq_1025}) and
(\ref{eq_925}), via Jensen's inequality.

\begin{lemma}  Let $d$ be a positive integer, $0 \leq \eps <
1$ and  $\mu \in \bM(S^{d-1})$. Suppose that we are
given a ``random probability measure'' on $S^{d-1}$. That is, let
$\lambda$ be a probability measure on a measurable space $\Omega$,
and suppose that with any $\alpha \in \Omega$ we associate a
probability measure $\mu_{\alpha} \in \bM(S^{d-1})$ such that the
map $ \Omega \ni \alpha \mapsto \mu_{\alpha} \in \bM(S^{d-1}) $ is
measurable. Assume that $$ d_{TV} \left( \mu, \int_{\Omega}
\mu_{\alpha} d \lambda(\alpha) \right) \leq \eps.
$$ Then,
$$ W_1(\mu, \sigma_{d-1}) \leq \int_{\Omega} W_1(\mu_{\alpha}, \sigma_{d-1}) d \lambda(\alpha) + \pi \eps
\leq \sup_{\alpha \in \Omega}  W_1(\mu_{\alpha}, \sigma_{d-1}) + 4 \eps. $$
\label{lem_sum}
\end{lemma}

Recall that $\mu|_{X}$ stands for the conditioning of $\mu$ to $X$.

\begin{lemma}
Suppose that $\mu$ and $\nu$ are finite measures on a measurable
space $\Omega$ and let $\eps > 0$. Let $X \subseteq \Omega$ be such
that $\nu(X) > \eps$. Suppose that
$$ |\mu(A) - \nu(A)| \leq \eps \ \ \ \ \ \ \text{for all} \ A \subseteq X. $$
 Then $ d_{TV}(\mu|_{X}, \nu|_{X}) \leq 2 \eps / \nu(X). $\label{useful}
\end{lemma}

\emph{Proof:} For any $A \subseteq X$,
$$ \big| \nu|_{X}(A) - \mu|_{X} (A) \big| =
\left| \frac{\nu(A) - \mu(A)}{\nu(X)} + \frac{\mu(A)}{\mu(X)} \cdot
\frac{\mu(X) - \nu(X)}{\nu(X)} \right| \leq \frac{\eps}{\nu(X)} +
\frac{\eps}{\nu(X)} = \frac{2 \eps}{\nu(X)}, $$ since $\mu(A) \leq
\mu(X)$.  \hfill $\square$

\medskip
Next we describe a few simple properties of $\eps$-radial measures.

\begin{lemma} Let $d$ be a positive integer and $0 < \eps < 1/2$.
Let $\mu, \nu$ be Borel probability measures on $\RR^d$.
Additionally, assume that we are given a ``random probability
measure'' on $\RR^d$. That is, let $\lambda$ be a probability
measure on a measurable space $\Omega$. Suppose that with any
$\alpha \in \Omega$ we associate a measure $\mu_{\alpha} \in
\bM(\RR^{d})$ such that the map $ \Omega \ni \alpha \mapsto
\mu_{\alpha} \in \bM(\RR^{d}) $ is measurable.
\begin{enumerate}
 \item[(a)] Suppose that $\mu$ is $\eps$-radial and that $d_{TV}(\mu, \nu) \leq \eps^2$.
Then $\nu$ is $5 \eps$-radial.
\item[(b)] Suppose that $\mu_{\alpha}$ is $\eps$-radial
for any $\alpha \in  \Omega$.
Assume that $\mu = \int_{\Omega} \mu_{\alpha} d \lambda(\alpha)$. Then $\mu$ is $4\sqrt{\eps}$-radial.
\item[(c)] Suppose that $A \subseteq \Omega$ satisfies $\lambda(A) \geq 1 - \eps$,
and $\mu_{\alpha}$ is $\eps$-radial for any $\alpha \in A$. Assume
that $\mu = \int_{\Omega} \mu_{\alpha} d \lambda(\alpha)$. Then
$\mu$ is $20\sqrt{\eps}$-radial.
\end{enumerate}
\label{lem_554}
\end{lemma}

\emph{Proof:} Begin with (a). Let $J \subset (0, \infty)$ be an
interval with $\nu(S(J)) \geq 2 \eps$, where $S(J) = \{ x \in \RR^d
; |x| \in J \}$ is a spherical shell. Denote $\nu_J = \nu|_{S(J)}$
and $\mu_J = \mu|_{S(J)}$. Since $d_{TV}(\mu, \nu) \leq \eps^2$, we
may apply Lemma \ref{useful} with $\eps^2$ and $X = S(J)$. We
conclude that $ d_{TV}(\mu_J, \nu_J) \leq 2 (\eps^2) / 2 \eps =
\eps. $ Consequently,
\begin{equation}
d_{TV}(\cR_*(\mu_J), \cR_*(\nu_J)) \leq \eps.
\label{eq_1005}
\end{equation}
Since $\mu$ is $\eps$-radial and $\mu(S(J)) \geq 2 \eps - \eps^2 \geq \eps$, then $W_1(\cR_*(\mu_J), \sigma_{d-1}) \leq \eps$ according to Definition \ref{def_934}. From (\ref{eq_925}), (\ref{eq_1005}) and the triangle inequality, $W_1(\cR_*(\nu_J), \sigma_{d-1}) \leq (\pi + 1) \eps \leq 5 \eps$. This completes the proof of (a).

\medskip We move to the proof of (b).
Let $J \subset (0, \infty)$ be an interval with
 $\mu(S(J)) \geq 4 \sqrt{\eps}$.
Let $X = \{ \alpha \in \Omega ; \mu_{\alpha}(S(J)) \geq \eps \}$.
Denote $\nu = \int_{X} \mu_{\alpha} d \lambda (\alpha)$, a finite
Borel measure on $\RR^n$. Then for any $A \subseteq S(J)$,
\begin{equation} \left| \mu(A) - \nu(A) \right| = \int_{\Omega \setminus X}\mu_{\alpha}(A) d \lambda
(\alpha) \leq \int_{\Omega \setminus X} \mu_{\alpha}(S(J)) d
\lambda(\alpha) \leq \eps \lambda(\Omega \setminus X) \leq \eps.
\label{eq_1234}
\end{equation}
Denote $\mu_J = \mu|_{S(J)}$ and $\nu_J = \nu|_{S(J)}$.
From (\ref{eq_1234}) and Lemma \ref{useful},
\begin{equation}
 d_{TV}(\mu_J, \nu_J) \leq 2 \eps / (4 \sqrt{\eps}) = \sqrt{\eps} / 2. \label{eq_428}
\end{equation}
Note that $\nu_J = \int_{X} \mu_{\alpha}|_{S(J)} d
\lambda^{\prime}(\alpha)$ where $\lambda^{\prime}$ is a certain
probability measure on $X$. Since $\mu_{\alpha}$ is $\eps$-radial
and $\mu_{\alpha}(S(J)) \geq \eps$ for $\alpha \in X$, then from
Definition \ref{def_934},
\begin{equation}
 W_1 \left( \cR_*(\mu_{\alpha}|_{S(J)}), \sigma_{d-1} \right) \leq \eps \ \ \ \ \ \ \text{for}
\ \alpha \in X. \label{eq_429}
\end{equation}
We have $\cR_*(\nu_J) = \int_{X} \cR_*(\mu_{\alpha}|_{S(J)}) d
\lambda^{\prime}(\alpha)$. Thus, (\ref{eq_429}) and Lemma
\ref{lem_sum} yield  that $W_1 \left(\cR_*(\nu_J), \sigma_{d-1} \right)
\leq \eps$. Combining the last inequality with (\ref{eq_925}) and
(\ref{eq_428}), we see that
$$ W_1 \left( \cR_*(\mu_J), \sigma_{d-1} \right) \leq
4 d_{TV} \left( \cR_*(\mu_J), \cR_*(\nu_J) \right)
+ W_1 \left(\cR_*(\nu_J), \sigma_{d-1} \right)
\leq 2 \sqrt{\eps} + \eps \leq 4 \sqrt{\eps}. $$
Since $\mu_J = \mu|_{S(J)}$, and $J \subset (0, \infty)$ is an arbitrary
interval with $\mu(S(J)) \geq 4 \sqrt{\eps}$, the assertion (b) is proven.

\medskip To prove (c), denote $\nu = \int_{A} \mu_\alpha d \lambda |_{A}(\alpha)$,
a probability measure on $\RR^d$. Then $\nu$ is
$4\sqrt{\eps}$-radial, according to (b). Furthermore, clearly
$d_{TV}(\lambda|_{A}, \lambda) = 1 - \lambda(A)  \leq  \eps$, and
hence $d_{TV}(\mu, \nu) \leq  \eps \leq (4 \sqrt{\eps})^2$.
According to part (a), the measure $\mu$ is $20 \sqrt{\eps}$-radial,
and (c) is proven. \hfill $\square$

\medskip
Probability measures are the protagonists of this text. Some of our
constructions of probability measures are probabilistic in nature.
To avoid confusion, we try to distinguish sharply between the
measures themselves, and the randomness used in their construction.
Whenever we have objects that are declared random (for instance,
random vectors in $S^{d-1}$), all statements containing probability
estimates or using the symbol $\PP$ refer to these random objects
and only to them.

\medskip The crude bound in the following lemma is
certainly a standard application of the so-called ``empirical
distribution method'' (see, e.g., \cite{blm_acta}). We were not able
to find it in the literature, so a proof is provided. Recall that
$\delta_x$ stands for the Dirac unit mass at the point $x$.

\begin{lemma} Let $d, N$ be positive integers,
and let $X_1,\ldots,X_N$ be independent random vectors, distributed
uniformly on $S^{d-1}$. Denote $ \mu = N^{-1} \sum_{i=1}^N \delta_{X_i}$. Then,
with probability greater than $1 - C\exp(-c \sqrt{N})$ of selecting $X_1,\ldots,X_N$,
$$ W_{1}(\mu, \sigma_{d-1}) \leq \frac{C}{N^{c/d}} $$
where $C > 1$ and $0 < c < 1$ are universal constants.
\label{lem_456}
\end{lemma}

\emph{Proof:} Denote by $\cF$ the class of all $1$-Lipschitz
functions $\vphi: S^{d-1} \rightarrow \RR$ such that  $\int \vphi d
\sigma_{d-1} = 0$. Note that $\sup | \vphi | \leq \pi$ for any
$\vphi \in \cF$. According to (\ref{eq_1455}),
\begin{equation}
W_1(\mu, \sigma_{d-1}) = \sup_{\vphi \in \cF} \int_{S^{d-1}} \vphi d
\mu. \label{eq_1044}
\end{equation}
 Let $\eps > 0$ be a parameter to be specified later on.
A subset $\cN \subset S^{d-1}$
is an $\eps$-net
if for any $x \in S^{d-1}$ there exists $y \in \cN$ with $\rho(x,y) \leq \eps$.
 Let $\cN$ be an $\eps$-net of cardinality
$\#(\cN) \leq (C / \eps)^d$ (see, e.g., \cite[Lemma 4.16]{Pisier}).
For $\vphi \in \cF$ denote
$$ \tilde{\vphi}(x) = \min_{y \in \cN} \left(  \eps \lceil \vphi(y) / \eps \rceil + \rho(x,y) \right), $$
where $\lceil a \rceil$ is the minimal integer that is not smaller
than $a$. Then $\tilde{\vphi}$ is a $1$-Lipschitz function, as a
minimum of
  $1$-Lipschitz functions. It is
easily verified that $\vphi \leq \tilde{\vphi} \leq \vphi + 3 \eps$.
Denote $\vphi^{\circ}(x) = \tilde{\vphi}(x) - \int \tilde{\vphi}(y)
d \sigma_{d-1}(y)$. Then $\vphi^{\circ} \in \cF$ for any $\vphi \in
\cF$, and $\sup| \vphi^{\circ} - \vphi | \leq 3 \eps$. Hence,
\begin{equation}
 W_1(\mu, \sigma_{d-1}) = \sup_{\vphi \in \cF}
\int_{S^{d-1}} \vphi d \mu \leq  3 \eps + \sup_{\vphi \in \cF}
\int_{S^{d-1}} \vphi^{\circ} d \mu = 3 \eps + \sup_{\vphi \in \cF}
\frac{1}{N} \sum_{i=1}^N \vphi^{\circ}(X_i). \label{eq_232}
\end{equation}
Denote $\tilde{\cF} = \{\tilde{\vphi} ; \vphi \in \cF\}$ and
$\cF^{\circ} = \{  \vphi^{\circ} ; \vphi \in \cF \}$. These sets are
finite. In fact, as each $\vphi \in \tilde{\cF}$ is determined by
the restriction $\vphi|_{\cN}$, we have
\begin{equation}
\#( \cF^{\circ}) \leq \#(\tilde{\cF}) \leq \left( \frac{2 \pi}{\eps}
+ 1 \right)^{\#(\cN)} \leq \exp \left( (C / \eps)^{2d} \right).
\label{eq_213}
 \end{equation}
Fix $\vphi^{\circ} \in \cF^{\circ}$. Then $\vphi^{\circ}$ is a
$1$-Lipschitz function on the sphere $S^{d-1}$ with $\int
\vphi^{\circ} d \sigma_{d-1} = 0$. According to L\'evy's lemma (see
Milman and Schechtman \cite[Section 2 and Appendix V]{MS}), for any
$i=1,\ldots,N$,
$$ \PP \left(
|\vphi^{\circ}(X_i)|  \geq t \right) \leq C \exp \left(-c t^2 d
\right) \ \ \ \ \ \forall t \geq 0, $$ where $\PP$ refers, of
course, to the probability of choosing the random vectors
$X_1,\ldots,X_N$. From Bernstein's inequality (see, e.g., Bourgain,
Lindenstrauss and Milman \cite[Proposition 1]{BLM}),
\begin{equation}
 \PP \left( \left| \frac{1}{N} \sum_{i=1}^N \vphi^{\circ}(X_i) \right|  \geq t \right) \leq C^{\prime} \exp
\left(-c^{\prime} t^2 N d \right) \ \ \ \ \ \forall t \geq 0.
\label{eq_229}
\end{equation}
Set $t = \eps$ in (\ref{eq_229}). From (\ref{eq_213})
and (\ref{eq_229}),
\begin{equation}
 \PP \left( \sup_{\vphi^{\circ} \in \cF^{\circ}}
\left| \frac{1}{N} \sum_{i=1}^N \vphi^{\circ}(X_i) \right|  \geq
\eps \right) \leq C^{\prime} \exp \left( (C / \eps)^{2d} -c^{\prime}
\eps^2 N d \right).  \label{eq_248}
\end{equation}
We now select $\eps = C N^{-1 / (2d+2)}$, for a sufficiently large
universal constant $C > 0$.  Substitute the value of $\eps$ in
(\ref{eq_248}) and apply (\ref{eq_232}), to deduce that
$$ W_1(\mu, \sigma_{d-1}) \leq 3 \eps + \sup_{\vphi^{\circ} \in \cF^{\circ}}
\frac{1}{N} \sum_{i=1}^N \vphi^{\circ}(X_i) \leq 4 \eps \leq
C^{\prime} N^{-1/(2d + 2)},
$$ with probability greater than $1 - C^{\prime} \exp(-c^{\prime}
N^{d / (d+1)})$ of selecting $X_1,\ldots,X_N$. \hfill $\square$

\medskip
\emph{Remark.} The discrepancy of $\mu \in \bM(S^{d-1})$ is defined as
$$ D(\mu) = \sup_{B} \left| \mu(B) - \sigma_{d-1}(B) \right| $$
where the supremum runs over all geodesic balls $B \subseteq
S^{d-1}$. A result analogous to Lemma \ref{lem_456} for discrepancy
appears in Beck and Chen \cite[Section 7.4]{beck_chen}. It is
possible to adapt our technique to suit the discrepancy metric, and
some of its variants, in place of $W_1$. In fact, the only
properties of the metric $W_1$ that are used in our proof are Lemma
\ref{lem_456} and (\ref{eq_1025}) and (\ref{eq_925}) above. Our
method, of course, works for the $W_p$ metrics as long as $1 \leq p
< \infty$, but it does not seem to apply for the $W_{\infty}$
metric. The $W_{\infty}$ metric induces a topology that is much
stronger than weak convergence, and it is not even weaker than
convergence in norm.

\section{Isotropic Gaussians}
\label{sec3}

A centered Gaussian random vector in $\RR^d$ is a random vector
whose density is proportional to $x \mapsto \exp(-M x \cdot x)$ for
a positive definite matrix $M$. A centered Gaussian random vector is
said to be isotropic if $M$ is a scalar matrix. It is called standard if $M =
Id /2$, where $Id$ is the identity matrix. Recall that $\cR$ stands
for radial projection.

\begin{lemma} Let $d, N$ be positive integers and let $Z_1,\ldots,Z_N$  be independent
isotropic Gaussian random vectors in $\RR^d$. Denote $\mu =
\frac{1}{N} \sum_{i=1}^N \delta_{Z_i}$.

\medskip Then, with probability greater than $1 - C \exp(-c N^{1/4})$
of selecting the $Z_i$'s,
the measure $\mu$ is $\delta$-radial, for $\delta = C N^{-c / d}$.
Here, $C > 1$ and $0 < c < 1$ are universal constants.
\label{lem_247}
\end{lemma}

\emph{Proof:}  Set $\eps = 5 / N^{1/4}$.
We may assume that $\eps \leq 1/ 10$, as otherwise
the conclusion of the lemma is obvious for a suitable
choice of universal constants $c, C > 0$.
 The central observation is that the radii $|Z_1|,\ldots,|Z_N|$
are independent of the angular parts $\cR(Z_1),\ldots, \cR(Z_N)$,
and that the random vectors $\cR(Z_1),\ldots, \cR(Z_N)$ are independent
random vectors that are distributed uniformly on $S^{d-1}$.

\medskip With probability one, none of the $|Z_i|$'s are zero,
and there are no $i \neq j$ with $|Z_i| = |Z_j|$.
We condition on the values $|Z_1|,\ldots, |Z_N|$, which are assumed
to be distinct and non-zero.
For an interval $J \subset (0, \infty)$  write
$$ Z(J) = \{ i ;  |Z_i| \in J \} \ \ \ \ \text{and} \ \ \ \ w(J) = \#(Z(J)). $$
Denote $k = \lceil 1 / \eps^2 \rceil$,
and let $J_1,\ldots,J_k \subset (0, \infty)$ be disjoint intervals
such that
$$ w(J_i) = \lfloor N i / k \rfloor -
\lfloor N (i-1) / k \rfloor \ \ \ \ \ \ \ \text{for} \ i=1,\ldots,k. $$
\begin{center}
\includegraphics[height=0.4in]{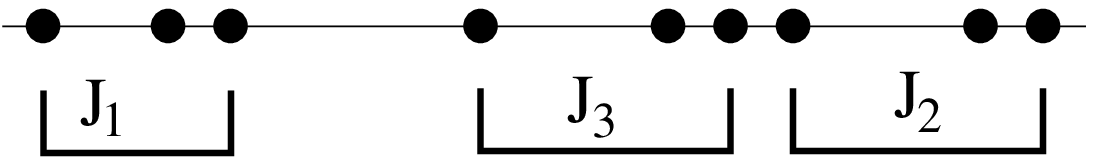}

{\small Figure 1}
\end{center}
Since $\eps^2 N \geq 2$,
then $\eps^2 N / 2 \leq w(J_i) \leq 2 \eps^2 N$
for any $i$.
For an interval $J \subset (0, \infty)$ with $w(J) \neq 0$,
denote $$ \mu_J = \frac{1}{w(J)} \sum_{j \in Z(J)} \delta_{\cR(Z_j)}. $$
Fix $i=1,\ldots,k$. We abbreviate $\mu_i = \mu_{J_i}$.
Since $\{ \cR(Z_j) \}_{j \in Z(J_i)}$ is a collection of $w(J_i)$
independent random vectors,  distributed uniformly on the sphere,
then Lemma \ref{lem_456} applies, and yields,
$$ \PP \left(  W_{1}(\mu_{i}, \sigma_{d-1}) \leq \frac{C}{w(J_i)^{c/d}} \right)
\geq 1 - C\exp(-c \sqrt{w(J_i)}). $$
We now let $i$ vary. Since $w(J_i)$ has the order of magnitude of $\eps^2 N$, then,
\begin{equation}
 \PP \left(  \max_{i=1,\ldots,k} W_{1}(\mu_i, \sigma_{d-1}) \leq \frac{C}{(\eps^2 N)^{c/d}}
 \right)
\geq 1 - C k \exp(-c \eps \sqrt{ N}).
\label{eq_932_}
\end{equation}
Write $\cI$ for the collection of all non-empty intervals in $(0, \infty)$.
Fix an interval $J$ with $w(J) \geq N \eps$. Let $J_{i_1},\ldots,J_{i_{\ell}}$
be all the intervals among the $J_i$'s that are contained in $J$. Then
$J_{i_1} \cup \ldots \cup J_{i_{\ell}}$ covers all but
at most $4 \eps^2 N$ of the $|Z_i|$'s that belong to the interval $J$.
Therefore,
$$ d_{TV} \left( \, \mu_J ,  \sum_{j=1}^{\ell} \lambda_j \mu_{i_j} \right)
\leq \frac{4 \eps^2 N}{N \eps } = 4 \eps
$$ where $\lambda_1,\ldots,\lambda_{\ell}$ are appropriate
non-negative coefficients that add to one. According to Lemma
\ref{lem_sum}, $$ W_1 \left(\mu_J, \sigma_{d-1} \right)  \leq
\max_{i=1,\ldots,k} W_1(\mu_{i}, \sigma_{d-1}) + 20 \eps  \ \ \ \ \
\ \text{for all} \ J \in \cI \ \text{with} \ w(J) \geq N \eps. $$ We
thus conclude from (\ref{eq_932_}) that with probability at least $1
- C k \exp(-c \eps \sqrt{N})$,
\begin{equation}
W_1(\mu_J, \sigma_{d-1}) \leq \frac{C}{(\eps^2 N)^{c / d}} + 20 \eps
\ \ \ \ \ \ \text{for all} \ J \in \cI \ \text{with} \ w(J) \geq N
\eps. \label{eq_939}
\end{equation}
The latter probability bound is valid under the conditioning on  $|Z_1|,\ldots,|Z_N|$, and it holds for all possible values of $|Z_1|,\ldots,|Z_N|$, up to measure zero.
Hence the aforementioned probability bound for (\ref{eq_939})
 also holds with no conditioning at all. Recall that we write $S(J) = \{ x \in \RR^n ; |x| \in J \}$,
and note that $\mu_J = \cR_*(\mu|_{S(J)})$ and $w(J) = N \mu(S(J))$.
Since $\eps \sqrt{N} \geq N^{1/4}$, then (\ref{eq_939}) translates
as follows: With probability greater than $1 - C \exp(-c N^{1/4})$
of selecting $Z_1,\ldots,Z_N$,
$$ W_1(\cR_*(\mu|_{S(J)}), \sigma_{d-1}) \leq C N^{-c/d} + C \eps
\ \ \ \ \ \ \text{for all} \ J \in \cI \ \text{with} \ \mu(S(J)) \geq \eps . $$
This means that $\mu$ is $C(N^{-c/d} +  \eps)$-radial with probability greater than $1 - C \exp(-c N^{1/4})$. Since
$N^{-c/d} +  \eps \leq C^{\prime} N^{-c^{\prime}/d}$, the lemma is proven.
\hfill $\square$

\begin{lemma} Let $k$ be a positive integer.
For an invertible
$k \times k$ matrix $A$, write $\gamma_A$
for the probability measure on $\RR^k$
whose density is proportional to $x \mapsto \exp(-|A x|^2 / 2)$.
Then, for any $k \times k$ invertible matrices $A$ and $B$,
$$ d_{TV}(\gamma_A, \gamma_B) \leq C k \| B A^{-1} - Id \| $$
where $Id$ is the identity matrix, $\| \cdot \|$ is the operator
norm, and $C > 0$ is a universal constant.
\label{lem_626}
\end{lemma}

\emph{Proof:} Let $X$ be a standard Gaussian random vector in $\RR^k$.
Then $\gamma_A(U) = \PP \left( A^{-1} X \in U\right)$ for any measurable set $U \subseteq \RR^k$. Therefore,
$$ d_{TV}(\gamma_A, \gamma_B) = \sup_{U \subseteq \RR^K} \left| \PP(X \in U) -
\PP(A B^{-1} X \in U )  \right| = d_{TV}(\gamma_{ B A^{-1} },
\gamma_{Id}). $$ Denote $M =  B A^{-1}$, write $\gamma =
\gamma_{Id}$ and set $\eps = \| M - Id \| = \sup_{|x| = 1} |Mx -
x|$. We write $\vphi_M(x) = (\det M) (2 \pi)^{-d/2} \exp(-|Mx|^2/2)$
for the density of $\gamma_M$ and similarly $\vphi$ stands for the
density of $\gamma$. We may assume that $\eps < 1/2$, as otherwise
the conclusion of the lemma is trivial. Then $\left| |M x|^2 - |x|^2
\right| \leq 3 \eps |x|^2$ for any $x \in \RR^k$, and also $(1 + 2
\eps)^{-k} \leq \det M \leq (1 + \eps)^k$. Therefore,
\begin{eqnarray*}
 \lefteqn{\left| \vphi(x) - \vphi_M(x) \right|} \\ & = &
\vphi(x) \left| 1 - (\det M) \exp \left(\frac{|x|^2 - |Mx|^2}{2}
\right)  \right| \leq \vphi(x) \left[ (1 + 2 \eps)^{k} \exp (3 \eps
|x|^2) - 1 \right]
\end{eqnarray*}
for any $x \in \RR^k$. Consequently,
$$ d_{TV} \left( \gamma, \gamma_M \right) = \frac{1}{2} \int_{\RR^k}
|\vphi(x) - \vphi_M(x)| dx \leq (1 + 2 \eps)^{k} \int_{\RR^k} \exp
(3 \eps |x|^2 ) \vphi(x) dx - 1. $$ However,
$$ \int_{\RR^k} \exp (3 \eps |x|^2 )
\vphi(x) dx = (2 \pi)^{-d/2} \int_{\RR^k} \exp \left( -|\sqrt{1 - 6
\eps} x|^2 / 2 \right) dx = (1 - 6 \eps)^{-k/2}.
$$
We deduce that
$$ d_{TV} \left( \gamma, \gamma_M \right) \leq   (1 + 2 \eps)^{k}
(1 - 6 \eps)^{-k/2} - 1 \leq C k \eps, $$ under the legitimate
assumption that $\eps < c / k$ (otherwise, there is nothing to
prove). \hfill $\square$

\medskip Consider the map $\cE: (\RR^d)^N \rightarrow \bM(\RR^d)$ defined by
$$ (\RR^d)^N \ni (x_1,\ldots,x_N) \stackrel{\cE}{\mapsto} \frac{1}{N} \sum_{i=1}^N \delta_{x_i}. $$
 A Borel probability measure $\alpha$ on $(\RR^d)^N$ thus
induces the Borel probability measure $\cE_*{\alpha}$ on the space
$\bM(\RR^d)$.
The next lemma is a small perturbation of Lemma \ref{lem_247}.

\begin{lemma} Let $d, N$ be positive integers, $\eps > 0$.
Let $X_1,\ldots,X_N$
be independent, standard Gaussian random vectors in $\RR^d$. Let $(a_{ij})_{1 \leq j\leq i \leq N}$ be real numbers, with $a_{ii} \neq 0$ for all $i$,
such that
\begin{equation}
 |a_{ij}| \leq \eps |a_{ii}|  \ \ \ \ \ \ \ \ \text{for} \ \ j < i.
\label{eq_647}
\end{equation}
Denote $Y_i = \sum_{j \leq i} a_{ij} X_j$ and consider the probability
measure $\mu = N^{-1} \sum_{i=1}^N \delta_{Y_i}$.

\medskip Then, with probability greater than $1 - C \exp(-c N^{1/4}) - C N^2 d^2 \eps$ of selecting the random vectors $X_1,\ldots,X_N$,
the measure $\mu$ is $\delta$-radial for $\delta = C N^{-c/d}$.
Here, $C > 1$ and $0 < c < 1$ are universal constants.
 \label{lem_324} \end{lemma}

\emph{Proof:}  Denote $Z_i = a_{ii} X_i$.
 The $Z_i$'s are independent, isotropic, Gaussian random vectors in $\RR^d$.
Denote by $U \subseteq \bM(\RR^d)$ the collection  of all
$\delta$-radial probability measures, where $\delta = C N^{-c/d}$ is
the same  as in Lemma \ref{lem_247}. Let $\alpha$ be the probability
measure on $(\RR^d)^N$ that is the joint distribution of
$Z_1,\ldots,Z_N$. According to Lemma \ref{lem_247},
$$ (\cE_* \alpha)(U) \geq 1 - C \exp(-c N^{1/4}).
$$
Let $\beta$ be the probability measure on $(\RR^d)^N$ which is the
joint distribution of $Y_1,\ldots,Y_N$. To prove the lemma, we need
to show that $$ (\cE_*\beta)(U) \geq 1 - C \exp(-c N^{1/4}) -
C^{\prime} N^2 d^2 \eps. $$ This would follow, if we could prove
that
\begin{equation}
 d_{TV}\left( (\cE_*\alpha), (\cE_*\beta) \right) \leq d_{TV}(\alpha, \beta) \leq C^{\prime} N^2 d^2 \eps.
\label{eq_640}
\end{equation}
Let $k = d N$. Let $A$ be the $k \times k$ matrix that represents
the linear operator
$$ \RR^k = (\RR^d)^N \ni (x_1,\ldots,x_N) \mapsto (a_{11} x_1,\ldots,a_{NN} x_N) \in (\RR^d)^N = \RR^k. $$
That is, $A$ is a diagonal matrix, and the diagonal of $A$ contains each  number $a_{ii}$ exactly
$d$ times. Denoting $X = (X_1,\ldots,X_N) \in (\RR^d)^N = \RR^k$ and
$Z = (Z_1,\ldots,Z_N) \in (\RR^d)^N = \RR^k$, we see that  $Z = A X$.
Therefore, in the notation of Lemma \ref{lem_626}, we have $\alpha = \gamma_{A^{-1}}$.

\medskip
Similarly, let $B$ be the $k \times k$ matrix that corresponds to
the linear operator
$$ (x_i)_{i=1,\ldots,N} \mapsto (\sum_{j \leq i} a_{ij} x_j)_{i=1,\ldots,N} $$
where $x_1,\ldots,x_N$ are vectors in $\RR^d$. Denoting $Y = (Y_1,\ldots,Y_N) \in \RR^k$
we have  $Y  = B X$
and consequently $\beta = \gamma_{B^{-1}}$.
Condition (\ref{eq_647}) implies that the off-diagonal elements of $A^{-1} B $
do not exceed $\eps$ in absolute value. The diagonal elements of $A^{-1} B$
are all ones.  Hence $\| A^{-1} B - Id \| \leq k \eps$, and according to Lemma  \ref{lem_626},
$$ d_{TV}(\alpha, \beta) = d_{TV}(\gamma_{B^{-1}}, \gamma_{A^{-1}})
\leq C k \| A^{-1} B - Id \| \leq C k^2 \eps = C N^2 d^2 \eps. $$
Thus (\ref{eq_640}) holds and the proof is complete.  \hfill $\square$

\section{Orthogonal decompositions}
\label{sec4}

This section is concerned with probability measures on $\RR^n$ that
may be  decomposed as a mixture,  whose components are mostly
ensembles of approximately-orthogonal vectors. Later on, we will
apply a random projection, and use Lemma \ref{lem_324} in order to
show that the projection of most elements in the mixture is
typically $\eps$-radial for small $\eps > 0$.

\begin{definition} Let $\ell, n$ be positive integers
and let $\eps > 0$.
Let $v_1,\ldots,v_{\ell} \in \RR^n$ be non-zero vectors,
and consider $v = (v_1,\ldots, v_{\ell})$, an $\ell$-tuple of vectors.
We say that $v$ is  ``$\eps$-orthogonal'' if
there exist orthonormal vectors $w_1,\ldots,w_{\ell} \in \RR^n$
and real numbers $(a_{ij})_{1 \leq j \leq i \leq \ell}$ such that
$$ v_i = \sum_{j=1}^i a_{ij} w_j \ \ \ \ \ \ \text{for} \ i=1,\ldots,\ell $$
and $|a_{ij}| \leq \eps |a_{ii}|$ for $j < i$.
 \label{def_327}
\end{definition}
Note that $(v_1,\ldots,v_{\ell})$ is $\eps$-orthogonal if and only
if $( \cR(v_1),\ldots,\cR(v_{\ell}))$ is $\eps$-orthogonal. We write
$\cO_{\ell,\eps} \subset (\RR^n)^{\ell}$ for the collection of all
$\eps$-orthogonal $\ell$-tuples $v = (v_1,\ldots,v_{\ell}) \in
(\RR^n)^{\ell}$. For a subspace $E \subset \RR^n$ denote by $Proj_E$
the orthogonal projection operator onto $E$ in $\RR^n$.

\begin{lemma} Let $\ell, n$ be positive integers. Suppose
that $\mu$ is a Borel probability measure on the unit sphere $S^{n-1}$ such that for any unit vector $\theta \in S^{n-1}$,
$$ \int_{S^{n-1}} (x \cdot \theta)^2 d \mu(x) \leq \frac{1}{\ell^{50}}. $$
Let $X_1,\ldots,X_\ell$ be independent random vectors in $S^{n-1}$, all distributed according to $\mu$. Then, with positive probability,
 $(X_1,\ldots,X_\ell)$ is $(\ell^{-20}/2)$-orthogonal.
\label{lem_124}
\end{lemma}

\emph{Proof:} We may assume that $\ell \geq 2$, otherwise
the lemma is vacuously true.
 Write $W_1,\ldots,W_{\ell} \in \RR^{n}$
for the vectors obtained from $X_1,\ldots,X_{\ell}$ through the
Gram-Schmidt orthogonalization process. (If the $X_i$'s are not
linearly independent, then some of the $W_i$'s might be zero). For
$i \geq 2$ denote by $E_i$ the linear span of $X_1,\ldots,X_{i-1}$.
Then, for $i \geq 2$,
$$ \EE |Proj_{E_i}(X_i)|^2 = \EE \sum_{j=1}^{i-1} (X_i \cdot W_j)^2
=  \EE \sum_{j=1}^{i-1} \int_{S^{n-1}} \left( x \cdot W_j \right)^2 d\mu(x)
\leq \frac{i-1}{\ell^{50}} \leq \ell^{-49}, $$
as $X_i$ is independent of $W_1,\ldots,W_{i-1}$.
 By Chebyshev's inequality,
$$ \PP \left \{ \exists 2 \leq i \leq \ell ;  |Proj_{E_i}(X_i)| \geq \ell^{-21} / 2 \right \} \leq (\ell-1) \frac{\ell^{-49}}{\ell^{-42} / 4}
\leq 4 \ell^{-6} < 1. $$
Therefore, with positive probability, $|Proj_{E_i}(X_i)| < \ell^{-21} / 2$ for all $i \geq 2$.
In this event, the vectors $X_1,\ldots,X_\ell$ are linearly independent,
and $W_1,\ldots,W_{\ell}$ are orthonormal vectors. Furthermore, in this event
 $a_{ii} := \sqrt{1 - |Proj_{E_i}(X_i)|^2} \geq 1 - \ell^{-21} / 2$
while $ a_{ij} := X_i \cdot W_j$ satisfies
$$ |a_{ij}| \leq |Proj_{E_i}(X_i)| \leq \ell^{-21} / 2  \ \ \ \ \ \ \ \ \ \text{for} \  j < i. $$
Thus, with positive probability, $X_i = \sum_{j \leq i} a_{ij} W_j$ for all $i$, with $|a_{ij}| \leq (\ell^{-20} / 2) |a_{ii}|$
for $j < i$ and with $W_1,\ldots,W_\ell$ being orthonormal vectors. This completes the proof.
Note that the ``positive probability'' is in fact greater than $1 - \ell^{-5}$. \hfill $\square$

\medskip The next lemma is essentially a measure-theoretic
variant of a lemma going back to Bourgain, Lindenstrauss and Milman
\cite[Lemma 12]{BLM}, with the main difference being that the
logarithmic dependence is improved upon to a power law. For two
Borel measures $\mu$ and $\nu$ on a compact $K$ we say that $\mu
\leq \nu$ if
\begin{equation} \int_{K} \vphi d \mu \leq \int_K \vphi d \nu \ \ \ \ \text{for any continuous}  \ \vphi:K \rightarrow [0, \infty).
\label{eq_954} \end{equation} Recall that a point is {\it not} in
the support of a measure if and only if it has an open neighborhood
of measure zero. We abbreviate $\cO_{\ell} = \cO_{\ell,
\ell^{-20}}$. For $v = (v_1,\ldots,v_\ell) \in \cO_{\ell}$ denote $$
\mu_v = \frac{1}{\ell} \sum_{i=1}^{\ell} \delta_{v_i}, $$ a Borel
probability measure on $\RR^n$ (in the notation of the previous
section, $\mu_v = \cE ( v)$). When $K \subseteq \RR^n$, we write $
\cO_{\ell}(K) \subseteq \cO_{\ell}$ for the collection of all
$(v_1,\ldots,v_{\ell}) \in \cO_{\ell} $ with $v_i \in K$ for all
$i$. Then $\cO_{\ell}(K) = \cO_{\ell} \cap K^{\ell} \subseteq
(\RR^n)^{\ell}$, and it is straightforward to verify that
$\cO_{\ell}(K)$ is compact whenever $K \subset \RR^n$ is a compact
that does not contain the origin.

\begin{lemma} Let $\ell, n$ be positive integers and let $0 < \eps < 1/2$.
 Let $\mu$ be a Borel probability measure on $\RR^n$
with $\mu(\{ 0 \}) = 0$. Assume that
\begin{equation}
 \sup_{\theta \in S^{n-1}}  \int_{S^{n-1}} (x \cdot \theta)^2  d \cR_*\mu(x) < \frac{\eps}{\ell^{50}}.
 \label{eq_1001}
 \end{equation}
Then, there exists a Borel probability measure $\lambda$
on $\cO_{\ell}$ such that
\begin{equation}
 d_{TV} \left( \mu,  \int_{\cO_{\ell}} \mu_v \, d \lambda(v) \right) <  \eps.
\label{eq_1047}
\end{equation}
\label{lem_1120}
\end{lemma}

\emph{Proof:} Since $\mu(\{ 0 \}) = 0$ then for any $\delta > 0$ we
may find a large punctured ball $K = \{ x \in \RR^n ; r \leq |x|
\leq R \}$ with $0 < r < R$ such that $\mu(K) \geq 1 - \delta$. We
may assume that $\mu$ is supported on a compact set that does not
contain the origin: Otherwise, replace $\mu$ with $\mu|_K$ for a
large punctured ball $K \subset \RR^n$ with $\mu(K) \geq 1 - \delta$
as above, and observe that $d_{TV}(\mu, \mu|_K) \leq   \delta$, so
the effect of the replacement on the inequalities (\ref{eq_1001})
and (\ref{eq_1047}) is bounded by $ \delta$, which can be made
arbitrarily small. Write $K \subset \RR^n$ for the support of $\mu$,
a compact which does not contain the origin. Denote by $\cF$ the
collection of all Borel measures $\lambda$ supported on
$\cO_{\ell}(K)$ for which
\begin{equation}
 \int_{\cO_{\ell}} \mu_v d \lambda(v) \leq \mu. \label{eq_932}
 \end{equation}
The condition (\ref{eq_932}) defining $\cF$ is closed in the weak$^*$ topology.
Furthermore, $\lambda(\cO_{\ell}) \leq 1$ for all $\lambda \in \cF$
(use (\ref{eq_932})), and take $\vphi \equiv 1$ in the definition (\ref{eq_954})).
Hence $\cF$ is a weak$^*$ closed subset of the unit ball of  the Banach space of
signed finite Borel measures on the compact $\cO_{\ell}(K)$.
From the Banach-Alaoglu theorem, $\cF$ is compact in the weak$^*$ topology.
Therefore the continuous functional $\lambda \mapsto \lambda(\cO_\ell)$ attains
its maximum on $\cF$ at some $\lambda_0 \in \cF$.
Clearly $\lambda_0(\cO_{\ell}) \leq 1$.
To prove the lemma, it suffices to show that
\begin{equation}
 \lambda_0(\cO_{\ell}) > 1 - \eps.
\label{eq_926}
\end{equation}
Indeed, if (\ref{eq_926}) holds, then we may define a probability
measure $\lambda_1 = \lambda_0 + \tilde{\lambda}$, where
$\tilde{\lambda}$ is {\it any} Borel measure on $\cO_{\ell}$ of
total mass $1 - \lambda_0(\cO_{\ell})$. Then clearly $$ d_{TV}
\left( \mu, \int_{\cO_{\ell}} \mu_v d \lambda_1(v) \right) \leq
\tilde{\lambda}(\cO_{\ell}) < \eps, $$ and the lemma follows. We
thus focus on the proof of (\ref{eq_926}). Assume by contradiction
that (\ref{eq_926}) fails. Denote $\nu = \mu - \int_{\cO_{\ell}}
\mu_v d \lambda_0(v)$. Then $\nu$ is a non-negative Borel measure on
$K \subset \RR^n$, according to (\ref{eq_932}), and also $\nu \leq
\mu$. Moreover, $\nu(K) \geq \eps$, since we assume that
(\ref{eq_926}) fails. Denote $\tilde{\nu} = \nu / \nu(K)$, a
probability measure on $K \subset \RR^n$. Then $\tilde{\nu} \leq \nu
/ \eps$ and hence $\cR_*(\tilde{\nu}) \leq \cR_*(\nu) / \eps \leq \cR_*(\mu) / \eps$. For
any unit vector $\theta \in S^{n-1}$,
$$ \int_{S^{n-1}} (x \cdot \theta)^2 d \cR_* \tilde{\nu}(x)
\leq \eps^{-1} \int_{S^{n-1}} (x \cdot \theta)^2 d \cR_*\nu(x) \leq
\eps^{-1} \int_{\RR^n} (x \cdot \theta)^2 d \cR_* \mu(x) \leq
\frac{1}{\ell^{50}}, $$ from our assumption (\ref{eq_1001}). Lemma
\ref{lem_124} thus asserts the existence of
$\tilde{x}_1,\ldots,\tilde{x}_\ell \in S^{n-1}$ in the support of
$\cR_*(\tilde{\nu})$ such that $(\tilde{x}_1,\ldots,\tilde{x}_\ell)$
is $(\ell^{-20}/2)$-orthogonal. Consequently, there exist non-zero
vectors $x_1,\ldots,x_{\ell} \in \RR^n$ in the support of $\nu$ such
that $(x_1,\ldots,x_{\ell})$ is $(\ell^{-20}/2)$-orthogonal. Let
$U_1,\ldots, U_{\ell} \subset \RR^n$ be small open neighborhoods of
$x_1,\ldots,x_{\ell}$, respectively, such that
$$ (y_1,\ldots,y_{\ell}) \in  \cO_{\ell, \ell^{-20}} = \cO_{\ell} \ \ \ \ \ \  \text{for all} \ \ y_1 \in U_1,\ldots, y_{\ell} \in U_{\ell}. $$
The $U_i$'s are necessarily disjoint and $U_1 \times ... \times
U_{\ell} \subseteq \cO_{\ell}$. Denote $\eta =
\min_{i=1,\ldots,\ell} \nu(U_i)$. Then $\eta > 0$, since $U_i$ is an
open neighborhood of the point $x_i$, and the point $x_i$ belongs to
the support of $\nu$. We set $\nu_i = \nu|_{U_i}$, the conditioning
of $\nu$ to $U_i$. Then $\nu_i$ is a probability measure supported
on $K \subset \RR^n$, and
$$ \eta \nu_i \leq \nu(U_i) \nu_i \leq \nu = \mu - \int_{\cO_{\ell}} \mu_v d \lambda_0(v)
\ \ \ \ \ \text{for} \  i=1,\ldots,\ell. $$
Therefore, also
$$  \eta \int_{U_1 \times ... \times U_{\ell}} \left( \frac{1}{\ell} \sum_{i=1}^{\ell} \delta_{y_i}  \right) d\nu_1(y_1) ... d \nu_\ell(y_\ell) = \eta \frac{\sum_{i=1}^{\ell} \nu_i}{\ell} \leq \mu - \int_{\cO_{\ell}} \mu_v d \lambda_0(v). $$
Consequently, the non-negative
measure $\lambda = \lambda_0 + \eta (\nu_1 \times ... \times \nu_{\ell})$
on $\cO_{\ell}(K)$ satisfies (\ref{eq_932}). Hence $\lambda \in \cF$, but $\lambda(\cO_{\ell}) =
\lambda_0(\cO_{\ell}) + \eta > \lambda_0(\cO_{\ell})$, in contradiction to the maximality of $\lambda_0$.
We thus conclude that (\ref{eq_926}) must be true, and the lemma is proven.
\hfill $\square$

\medskip A $d \times n$
matrix $\Gamma$ will be called a ``standard Gaussian random matrix''
if the entries of $\Gamma$ are independent standard Gaussian random
variables (of mean zero and variance one). Suppose that $w_1,\ldots,
w_{\ell}$ are orthonormal vectors in $\RR^n$ and that $\Gamma$ is a $d \times n$
standard Gaussian random matrix. Observe that in this
case, $\Gamma(w_1),\ldots, \Gamma(w_\ell)$ are independent standard
 Gaussian random vectors in $\RR^d$.

\begin{lemma} Let $d \leq \ell \leq n$ be positive integers,
let $0 < \eps < 1$ and assume that
\begin{equation}
 \ell \geq (C / \eps)^{Cd}.
\label{eq_1909} \end{equation}
  Suppose that $\lambda$ is a Borel
probability measure on $\cO_{\ell}$, and denote $\mu =
\int_{\cO_{\ell}} \mu_v \, d \lambda(v)$. Let $\Gamma$ be a $d
\times n$ standard Gaussian  random matrix.

\medskip Then, with positive probability of selecting the random matrix $\Gamma$,
the measure $\Gamma_* \mu$
on $\RR^d$ is $\eps$-radial. Here, $C > 1$
is a universal constant. (In fact, this probability is at least $1 - \ell^{-8}$.)
\label{lem_220}
\end{lemma}

\emph{Proof:} Fix $v = (v_1,\ldots, v_{\ell}) \in \cO_{\ell}$.
Consider the  measure $\tilde{\mu}_v := \Gamma_*(\mu_v)$ on $\RR^d$.
Then $\tilde{\mu}_v = \Gamma_*(\frac{1}{\ell} \sum_{i=1}^{\ell}
\delta_{v_i} ) = \frac{1}{\ell} \sum_{i=1}^{\ell}
\delta_{\Gamma(v_i)}$. Let $E(v)$ be the following event:
\begin{itemize}
 \item The measure $\tilde{\mu}_v$ is $C \ell^{-c/d}$-radial,
where $C > 1$ and $0 < c < 1$ are the universal constants from Lemma \ref{lem_324}.
\end{itemize}
Let us emphasize that for any $v \in \cO_{\ell}$, the event $E(v)$
might either hold or not, depending on the Gaussian random matrix $\Gamma$.
We are going to apply Lemma \ref{lem_324}. Since $v \in \cO_{\ell} = \cO_{\ell, \ell^{-20}}$,
then there exist orthonormal vectors $w_1,\ldots,w_{\ell} \in \RR^n$
and numbers $a_{ij}$ such that $v_i = \sum_{j \leq i} a_{ij} w_j$
and $|a_{ij}| \leq \ell^{-20} |a_{ii}|$ for $j < i$, with $a_{i i} \neq 0$
for all $i$. Denote $X_i = \Gamma(w_i)$
and $Y_i = \sum_{j \leq i} a_{ij} X_j$.
Then $X_1,\ldots,X_\ell$ are independent standard Gaussian random vectors
in $\RR^d$, and $\tilde{\mu}_v = \ell^{-1} \sum_{i=1}^\ell \delta_{Y_i}$.
We may thus apply Lemma \ref{lem_324} (with $N = \ell$ and $\eps = \ell^{-20}$)
 and conclude that for any $v \in \cO_{\ell}$,
\begin{equation}
 \PP(E(v)) \geq 1 - C \exp(-c \ell^{1/4}) - C \ell^2 d^2 \ell^{-20} \geq 1 - C^{\prime} \ell^{-16}. \label{event}
\end{equation}
Let $\cF \subseteq \cO_{\ell}$ be the collection of all $v \in \cO_{\ell}$
for which the event $E(v)$ holds. Then $\cF$ is a random
subset of $\cO_{\ell}$ (depending on the random matrix $\Gamma$).
According to (\ref{event}),
$$ \EE \lambda(\cF) = \EE \int_{\cO_{\ell}} 1_{\cF}(v) d \lambda(v) =
 \int_{\cO_{\ell}} \EE 1_{\cF}(v) d \lambda(v) =
\int_{\cO_{\ell}} \PP(E(v)) d \lambda(v) \geq 1 - C^{\prime} \ell^{-16}
 $$
where $1_{\cF}$ is the characteristic function of $\cF$.
Therefore, by Chebyshev's inequality,
$$ \PP \left( \lambda(\cF) \leq 1 - 2 C^{\prime} \ell^{-8}  \right)
\leq \frac{\EE \left(1 - \lambda(\cF) \right)}{2 C^{\prime}
\ell^{-8} } \leq  \ell^{-8} / 2 < 1. $$ We may assume that
$C^{\prime} \ell^{-8} \leq 1/2$, thanks to (\ref{eq_1909}). We
showed that with positive probability $\lambda(\cF) \geq 1 - 2
C^{\prime} \ell^{-8}$.  Recall that $\tilde{\mu}_v =
\Gamma_*(\mu_v)$ is $C \ell^{-c/d}$-radial for any $v \in \cF$.
Hence, according to Lemma \ref{lem_554}(c), with positive
probability of selecting the Gaussian random matrix $\Gamma$, the
measure
$$ \int_{\cO_\ell} \Gamma_*(\mu_v) d \lambda(v) = \Gamma_* \left(
\int_{\cO_\ell} \mu_v d \lambda(v) \right) = \Gamma_*(\mu) $$
is $C^{\prime} \ell^{-c^{\prime} / d}$-radial on $\RR^d$.
  \hfill $\square$

\medskip The Grassmannian $G_{n, k}$
of all $k$-dimensional subspaces in $\RR^n$ carries a unique
rotationally invariant probability measure, which will be referred
to as the uniform probability measure on $G_{n,k}$. When $\Gamma$ is
a $d \times n$ standard Gaussian  random matrix, the kernel of
$\Gamma$ is a random $(n-d)$-dimensional subspace, that is
distributed uniformly in the Grassmannian $G_{n, n-d}$. For a
subspace $E \subseteq \RR^n$ we write $E^{\perp} = \{ x \in \RR^n ;
\forall y \in E, x \cdot y = 0 \}$ for its orthogonal complement.

\begin{lemma} Let $0 \leq k \leq n-1$ be integers,
 and let $\mu$ be a Borel probability measure on $\RR^n$ with $\mu ( \{ 0 \} ) = 0$.
 Suppose that
 $E$ is a random $k$-dimensional subspace, distributed uniformly in $G_{n,k}$.
Then $\mu(E) = 0$ with probability one of selecting $E$.
\label{obvious}
\end{lemma}

\emph{Proof:} By induction on $k$. The case $k=0$ holds trivially.
Suppose now that $k \geq 1$, let $n$ be such that $k \leq n-1$, and
let $\mu$ be a Borel probability measure on $\RR^n$ with $\mu ( \{ 0
\} ) = 0$. Since $\mu ( \{ 0 \} ) = 0 $, there are at most countably
many one-dimensional subspaces $\ell \subset \RR^n$ with $\mu(\ell)
> 0$. Let $\ell$ be a random one-dimensional subspace, distributed
uniformly in $G_{n,1}$.  Then with probability one, $\mu(\ell) = 0$.
Denote $\nu = (Proj_{\ell^{\perp}})_* \mu$, a measure supported on
an $(n-1)$-dimensional subspace, with $\nu(\{ 0 \}) = 0$. Let $F$ be
a random $(k-1)$-dimensional subspace in $\ell^{\perp}$, distributed
uniformly. By the induction hypothesis, $\nu(F) = 0$ with
probability one. Denoting $E = Proj_{\ell^{\perp}}^{-1}(F)$, we see
that $\mu(E) = \nu(F) = 0$ with probability one. From the uniqueness
of the Haar measure, $E$ is distributed uniformly in $G_{n,k}$, and
the lemma follows.
 \hfill $\square$

\begin{corollary}
Let $1 \leq d \leq n$ be integers and let $0 < \eps < 1/2$.
 Let $\mu$ be a Borel probability measure on $\RR^n$
with $\mu(\{ 0 \}) = 0$. Assume that
\begin{equation}
 \sup_{\theta \in S^{n-1}} \int_{S^{n-1}} (x \cdot \theta)^2  d \cR_*\mu(x) \leq (\tilde{c} \eps)^{\tilde{C} d}.
 \label{eq_1142}
 \end{equation}
 Let $\Gamma$ be a $d \times n$ standard Gaussian random matrix.
Then, with positive probability of selecting $\Gamma$, the measure $\Gamma_* \mu $
on $\RR^d$ is $\eps$-radial proper. Here, $0 < \tilde{c} < 1$ and $\tilde{C} > 1$
are universal constants. (In fact, we have a lower bound of $1 - (\tilde{c} \eps)^{\tilde{C} d / 10}$ for the aforementioned probability.)
\label{cor_216}
\end{corollary}

\emph{Proof:} Throughout this proof, we write $C$ for the universal
constant from Lemma \ref{lem_220}. We define $\tilde{c} = (10
C)^{-1}$ and $\tilde{C} = 100 C$. It is elementary to verify that
with this choice of universal constants, there exists an integer
$\ell$ such that
$$ \ell \geq (5 C / \eps)^{C d} \ \ \ \ \ \text{and} \ \ \ \ \ (\tilde{c} \eps)^{\tilde{C} d} \leq \frac{\eps^2}{50 \ell^{50}}. $$
Note that the left-hand side of (\ref{eq_1142}) is at least $1/n$. Indeed,
\begin{eqnarray*}
 \lefteqn{\sup_{\theta \in S^{n-1}} \int_{S^{n-1}} (x \cdot \theta)^2  d \cR_*\mu(x) } \\
& \geq & \int_{S^{n-1}} \int_{S^{n-1}} (x \cdot \theta)^2  d \cR_*\mu(x) d \sigma_{d-1}(\theta)
= \int_{S^{n-1}} \frac{|x|^2}{n} d \cR_*(\mu)(x) = \frac{1}{n}.
\end{eqnarray*}
We conclude that $d < \ell \leq n$, and
$$
 \sup_{\theta \in S^{n-1}} \int_{S^{n-1}} (x \cdot \theta)^2  d \cR_* \mu(x) \leq
 \frac{\eps^2}{50 \ell^{50}}.
 $$
According to Lemma \ref{lem_1120},
there exists a Borel probability measure $\lambda$
on $\cO_{\ell}$ such that
\begin{equation}
 d_{TV} \left( \mu,  \int_{\cO_{\ell}} \mu_v \, d \lambda(v) \right) \leq  \eps^2 / 25.
\label{eq_218}
\end{equation}
Denote $\nu = \int_{\cO_{\ell}} \mu_v \, d \lambda(v)$.
From Lemma \ref{lem_220} the measure $\Gamma_*(\nu)$ is $(\eps/5)$-radial, with
probability at least $1-\ell^{-8}$ of selecting $\Gamma$,
because $\ell \geq (5 C / \eps)^{C d}$.
Additionally, $d_{TV}(\Gamma_* \mu, \Gamma_* \nu) \leq \eps^2 / 25$,
by (\ref{eq_218}). From Lemma \ref{lem_554}(a) we thus
learn that $\Gamma_*(\mu)$ is $\eps$-radial,
with positive probability of selecting $\Gamma$.
Moreover, $\Gamma_*(\mu)(\{ 0 \}) = \mu(\Gamma^{-1}(0)) = 0$
with probability one, according to Lemma \ref{obvious}.
 Hence, with positive probability,
$\Gamma_*(\mu)$ is $\eps$-radial proper. \hfill $\square$

\section{Selecting a position}
\label{sec5}

Our goal in this section is to find an appropriate invertible linear
transformation $T$ on $\RR^n$ such that $T_* \mu$ satisfies the
requirements of Corollary \ref{cor_216}. Our analysis is very much
related to the results of Barthe \cite{barthe_inventiones}, Carlen
and Cordero-Erausquin \cite{CC} and Carlen, Lieb and Loss
\cite{CLL}. For $x = (x_1,\ldots,x_n) \in \RR^n$ we write $x \otimes
x$ for the $n \times n$ matrix whose entries are $(x_i x_j)_{i,j =
1,\ldots,n}$. For a probability measure $\mu$ on the unit sphere
$S^{n-1}$ define
$$ M(\mu) = \int_{S^{n-1}} (x \otimes x) d \mu(x). $$
Then $M(\mu)$ is a positive semi-definite matrix of trace one.
Clearly, for any $\theta \in \RR^n$, we have $M(\mu) \theta \cdot
\theta = \int_{S^{n-1}} (x \cdot \theta)^2 d \mu(x)$. More
generally, for any  subspace $E \subseteq \RR^n$,
$$ \int_{S^{n-1}} |Proj_E(x)|^2 d \mu(x) = Tr(  Proj_E M(\mu) ) = Tr(  M(\mu) Proj_E), $$
the trace of the matrix $M(\mu) Proj_E $. A Borel probability
measure $\mu$ on $S^{n-1}$ is called isotropic if $M(\mu) = Id / n$,
where $Id$ is the identity matrix. Observe that when $\mu$ is
isotropic, for any subspace $E \subseteq \RR^n$,
\begin{equation}
\int_{S^{n-1}} |Proj_E x|^2 d \mu(x) = \frac{\dim(E)}{n}.
\label{eq_2142}
\end{equation} In particular, $\mu(E) \leq \dim(E) / n$ and
hence an isotropic probability measure is necessarily decent in the sense
of Definition \ref{def_1254}. A Borel probability measure $\mu$ on $\RR^n$ with $\mu(\{0 \}) = 0$ is
called ``potentially isotropic'' if there exists an invertible
linear map $T$ on $\RR^n$ such that $(\cR \circ T)_* \mu$ is
isotropic.

\begin{lemma} Let $\mu$ be a Borel probability measure on $S^{n-1}$ such that \begin{equation} \mu(H) = 0
\label{eq_1209}
\end{equation} for any hyperplane $H \subset \RR^n$ through the origin.
Then $\mu$ is potentially isotropic. \label{lem_1117}
\end{lemma}

\emph{Proof:} Given an invertible linear map $T: \RR^n \rightarrow
\RR^n$ we abbreviate $M_{\mu}(T) = M( (\cR \circ T)_* \mu )$. Then
$M_{\mu}(T)$ is a positive semi-definite matrix of trace one, and by
the arithmetic-geometric means inequality, $\det M_{\mu}(T) \leq n^{-n}$.
Note that $M_{\mu}(T) = M_{\mu}(\lambda T)$ for any $\lambda > 0$.
Consider the supremum of the continuous functional
\begin{equation} T \mapsto \det M_{\mu}(T) \label{eq_1213}
\end{equation} over the space of all
invertible linear operators $T: \RR^n \rightarrow \RR^n$ of
Hilbert-Schmidt norm one.

\medskip We claim that the supremum is attained. Indeed,
let $T_1,T_2,\ldots$ be a maximizing sequence of matrices. By
passing to a subsequence if necessary, we may assume that $T_i
\rightarrow T$, for a certain matrix $T$ of Hilbert-Schmidt norm
one. We need to show that $T$ is invertible. Denote by $E$ the image
of $T$, a subspace of $\RR^n$. We need to show that $E = \RR^n$. For
any $x \in S^{n-1}$ which is not in the kernel of $T$, we have $T_i
x \rightarrow T x \in E \setminus \{ 0 \}$, hence
\begin{equation} |Proj_{E^{\perp}} (\cR \circ T_i)(x)| \stackrel{i \rightarrow
\infty}{\longrightarrow} 0. \label{eq_1207}
\end{equation}
The kernel of $T$ is at most $(n-1)$-dimensional, since
the Hilbert-Schmidt norm of $T$ is one. According to
(\ref{eq_1209}), the convergence in (\ref{eq_1207}) occurs
$\mu$-almost-everywhere in $x$. Therefore,
$$ Tr( M_{\mu}(T_i) Proj_{E^{\perp}} ) =
\int_{S^{n-1}} |Proj_{E^{\perp}} (\cR \circ T_i)(x)|^2 d \mu(x)
\stackrel{i \rightarrow \infty}{\longrightarrow} 0. $$ We conclude
that if $E \neq \RR^n$, then $\det M_{\mu}(T_i) \rightarrow 0$, in
contradiction to the maximizing property of the sequence $(T_i)_{i
\geq 1}$. Hence $E = \RR^n$ and $T$ is invertible. Thus the supremum
of the functional (\ref{eq_1213}) is attained for some invertible
matrix $T_0$ of Hilbert-Schmidt norm one. We will show that $(\cR
\circ T_0)_* \mu$ is isotropic. Without loss of generality we assume
that $T_0 = Id$ (otherwise, replace $\mu$ with $(\cR \circ T_0)_*
\mu$ and note that this replacement does not affect the validity of
the assumptions and the conclusions of the lemma).

\medskip
The matrix $M(\mu) = M_{\mu}(Id)$ is a positive semi-definite matrix
of trace one. It is non-singular, thanks to (\ref{eq_1209}), and
therefore $M(\mu)$ is in fact positive definite. Moreover, for any
function $u: S^{n-1} \rightarrow \RR$ which is positive $\mu$-almost-everywhere and
for any $\theta \in S^{n-1}$,
\begin{equation}
\int_{S^{n-1}} u(x) (x \cdot \theta)^2 d \mu(x) > 0.
\label{eq_1754}
\end{equation}
Assume by contradiction that $M(\mu)$ is not a scalar matrix. Denote
by $\lambda_1$ the largest eigenvalue of $M(\mu) = M_{\mu}(Id)$, and
let $E \subset \RR^n$ be the eigenspace corresponding to the
eigenvalue $\lambda_1$. Then $1 \leq \dim(E) \leq n-1$. For $0 \leq
\delta < 1$ consider the linear operator
$$ L_{\delta}(x) = x - \delta Proj_E(X) \ \ \ \ \ \ \ (x \in \RR^n). $$
Then $Proj_E L_{\delta} = (1 - \delta) Proj_E$ while
$Proj_{E^{\perp}} L_{\delta} = Proj_{E^{\perp}}$.
This means that $\cR \circ L_{\delta}$ strengthens the $E^{\perp}$-component
of a given point in $\RR^n$, at the expense of its $E$-component. More precisely, for any $x
\in S^{n-1}$ and $0 \leq \delta < 1$ there exists $\eps_x^\delta
\geq 0$ such that
$$ Proj_{E^{\perp}} (\cR(  L_{\delta} x )) = (1 + \eps_{x}^{\delta}) Proj_{E^{\perp}} (x) . $$
Moreover, when $ x \not \in E \cup E^{\perp}$ we have the inequality
 $\eps_x^{\delta} \geq \eps(x)\delta$ for some $\eps(x) > 0$ depending only on $x$.
Consequently, for any $0 \leq \delta < 1$ and a non-zero vector
$\theta \in E^{\perp}$,
\begin{eqnarray}
\nonumber
 \lefteqn{
 \int_{S^{n-1}} (\cR(L_{\delta} x) \cdot \theta)^2 d \mu (x)
 = \int_{S^{n-1}} (1 + \eps_x^{\delta})^2 (x \cdot \theta)^2 d \mu
 (x)}\\ & \geq & \int_{S^{n-1}} (x \cdot \theta)^2 d \mu
(x) + 2 \delta \int_{S^{n-1}} \eps(x) (x \cdot \theta)^2 d \mu(x).
\label{eq_17541}
 \end{eqnarray} The symmetric matrix $M_{\mu}(L_{\delta})$
is of trace one, and it depends smoothly on $\delta$. Denote $D = d
M_\mu(L_{\delta}) / d\delta|_{\delta = 0}$, a traceless symmetric
matrix. According to our assumption (\ref{eq_1209}), the condition
$x \not \in E \cup E^{\perp}$ holds $\mu$-almost-everywhere as $1
\leq \dim(E) \leq n-1$. Therefore $\eps(x) > 0$ for $\mu$-almost
every $x \in S^{n-1}$. From (\ref{eq_1754}) and (\ref{eq_17541}) we
learn that for any $0 \neq \theta \in E^{\perp}$,
\begin{equation}
D \theta \cdot \theta = \left. \frac{d}{d \delta} \left( \int_{S^{n-1}} (\cR(L_{\delta} x) \cdot \theta)^2 d \mu (x) \right) \right|_{\delta = 0} > 0. \label{eq_1733}
\end{equation}
Recall that $E$ is the eigenspace corresponding to the maximal
eigenvalue $\lambda_1$ of $M(\mu)$. Denote by $\lambda_2$ the
second-largest eigenvalue, which is still positive but is strictly
smaller than $\lambda_1$. Then $Proj_{E^{\perp}} M(\mu)^{-1} \geq
\lambda_2^{-1} Proj_{E^{\perp}}$ in the sense of symmetric matrices.
Using elementary linear algebra, we deduce from (\ref{eq_1733}) that
$$ Tr( Proj_{E^{\perp}} M(\mu)^{-1} D )  \geq \lambda_2^{-1} Tr(
Proj_{E^{\perp}} D ) > 0. $$ Since $Tr(D) = 0$ then $Tr(Proj_E D) =
- Tr(Proj_{E^{\perp}} D)$ and
\begin{eqnarray*}
\lefteqn{ \left. \frac{d \log \det M_{\mu}(L_{\delta})}{d \delta}
\right|_{\delta = 0} = Tr \left( M(\mu)^{-1} D \right) } \\ & = & Tr(
Proj_{E} M(\mu)^{-1} D ) + Tr( Proj_{E^{\perp}} M(\mu)^{-1} D )
 \\ & = &
\frac{Tr (Proj_E D)}{\lambda_1}  + Tr( Proj_{E^{\perp}} M(\mu)^{-1} D )
 \geq  \left( \frac{1}{\lambda_2} - \frac{1}{\lambda_1} \right) Tr( Proj_{E^{\perp}} D ) > 0,
\phantom{aaaaa}
\end{eqnarray*}
in contradiction to the maximality of $\det M(\mu)$. Hence our
assumption that $M(\mu)$ is {\it not} a scalar matrix was absurd.
Since $Tr(M(\mu)) = 1$ then $M(\mu) = Id / n$ and $\mu$ is
isotropic.
 \hfill $\square$

\medskip
For a subspace $E \subset \RR^n$ and $\delta
> 0$ we write
\begin{equation} \cN_{\delta}(E) = \left \{ r x ; | x | = 1, r \geq
0, d(x, E \cap S^{n-1}) \leq \delta \right \} \label{eq_2050}
\end{equation}  where
$d(x, A) = \inf_{y \in A} | x - y |$. Then $\cN_{\delta}(E)$ is the
projective $\delta$-neighborhood of $E$. We will need the following
auxiliary continuity lemma. It is the only time in this text where
the non-degeneracy conditions of Definition \ref{def_1254} are used.

\begin{lemma} Let $n \geq 1$ be an integer
and let $\mu$ be a probability measure on $\RR^n$ with $\mu (\{0 \})
= 0 $ such that
$$ \mu(E) < \dim(E) / n $$
for any subspace $E \subset \RR^n$ other than $\RR^n$ and $\{ 0 \}$.
Suppose there exists a sequence of potentially isotropic probability
measures on $S^{n-1}$ that converges to $\cR_* \mu$ in the weak$^*$
topology. Then  $\mu$ is potentially isotropic. \label{lem_2019}
\end{lemma}

\emph{Proof:} From the assumptions of the lemma, there exist Borel
probability measures $\mu_1,\mu_2,\ldots$ on $S^{n-1}$ and
invertible linear maps $T_1,T_2,\ldots$ for which the following
holds:
\begin{itemize}
\item $\mu_i \stackrel{i \rightarrow \infty}{\longrightarrow} \cR_* \mu$
in the weak$^*$ topology; and
\item $(\cR \circ T_i)_* \mu_i$ is
isotropic for all $i$.
\end{itemize}
Without loss of generality we may assume that the $T_i$'s are
positive definite operators of trace one: If not, we will replace
the operator $T_i$ by $r U_i T_i$, where $U_i$ is an orthogonal
transformation such that $U_i T_i$ is positive definite and $r^{-1}$
is the trace of $U_i T_i$. Such a replacement does not affect the
isotropicity of $(\cR \circ T_i)_* \mu_i$. Furthermore, replacing
$T_i \ (i=1,2,\ldots)$ with a subsequence, we may assume that $T_i
\rightarrow T$, where $T$ is a positive semi-definite matrix of
trace one.

\medskip We claim that $T$ is invertible. Assume by contradiction
that $T$ is singular. Denote by $E \subset \RR^n$ the kernel of $T$,
and set $k = \dim(E)$. Then $1 \leq k \leq n-1$ and hence $\mu(E) <
k / n$. Since $E = \cap_{\delta > 0} \cN_{\delta}(E)$, then there
exists $\delta > 0$ such that $$ \mu(\cN_{\delta}(E) ) < k/n. $$ The
set $\cN_{\delta}(E) \cap S^{n-1}$ is closed in $S^{n-1}$. Since
$\mu_i \rightarrow \cR_* \mu$ in the weak$^*$ topology then
\begin{equation}
 \limsup_{i \rightarrow \infty} \mu_i( \cN_{\delta}(E) ) \leq
\cR_*\mu(\cN_{\delta}(E) ) = \mu(\cN_{\delta}(E) ) < \frac{k}{n}.
\label{eq_2113}
\end{equation}
Recall that $T_i \rightarrow T$, that the $T_i$'s are self-adjoint,
and that $E$ is the kernel of $T$, hence $E^{\perp}$ is the image of
$T$. This entails, roughly speaking, that for any $x \not \in E$,
the sequence $T_i x$ is ``approaching $E^{\perp}$''. In more precise
terms, we conclude that for any $x \not \in \cN_{\delta}(E)$,
\begin{equation}
 \left| Proj_{E^{\perp}} \left( \frac{T_i x}{|T_i x|} \right) \right| \stackrel{i
\rightarrow \infty} \longrightarrow 1. \label{eq_2107}
\end{equation}
Moreover, the convergence in (\ref{eq_2107}) is uniform over $x \in
\RR^n \setminus \cN_{\delta}(E)$. Consequently, from (\ref{eq_2113})
and (\ref{eq_2107}),
\begin{equation}
 \liminf_{i \rightarrow \infty}
 \int_{S^{n-1} \setminus \cN_{\delta}(E)}
\left| Proj_{E^{\perp}} \left( \frac{T_i x}{|T_i x|} \right)
\right|^2 d \mu_i(x) = \liminf_{i \rightarrow \infty} \mu_i( S^{n-1}
\setminus \cN_{\delta}(E) ) > 1 - \frac{k}{n}. \label{eq_2146}
\end{equation}
Recall that $(\cR \circ T_i)_* \mu_i$ is isotropic. According to
(\ref{eq_2142}), $$
 \int_{S^{n-1}} \left| Proj_{E^{\perp}} \left( \frac{T_i x}{|T_i
x|} \right) \right|^2 d \mu_i(x) = \frac{\dim(E^{\perp})}{n} = 1 -
\frac{k}{n},
$$
in contradiction to (\ref{eq_2146}). Thus our assumption that $T$ is
singular was absurd, and $T$ is necessarily invertible.

\medskip Since $T_i \rightarrow T$ with $T$ being invertible, we  know
that for any $x \in S^{n-1}$,
$$ \frac{T_i x}{|T_i x|} \stackrel{i \rightarrow
\infty}{\longrightarrow} \frac{T x}{|T x|} $$ and the convergence is
uniform in $S^{n-1}$. Therefore, for any $\theta \in S^{n-1}$,
\begin{equation}
  \int_{S^{n-1}} \left( \frac{T_i x}{|T_i x|} \cdot \theta \right)^2 d
\mu_i(x) \stackrel{i \rightarrow \infty}{\longrightarrow}
\int_{S^{n-1}} \left( \frac{T x}{|T x|} \cdot \theta \right)^2 d
\cR_* \mu(x). \label{eq_2154}
\end{equation}
However, the left-hand side of (\ref{eq_2154}) is always $1/n$. We
see that $(\cR \circ T \circ \cR)_* \mu = (\cR \circ T)_* \mu$  is
isotropic, and therefore $\mu$ is potentially isotropic. \hfill
$\square$

\begin{corollary} Let $n$ be a positive integer
and let $\mu$ be a Borel probability measure on $\RR^n$ with $\mu
(\{0 \}) = 0 $ such that
\begin{equation}
\mu(E) < \dim(E) / n
\label{eq_1039} \end{equation}
for any subspace $E \subset \RR^n$ other than $\RR^n$ and $\{ 0 \}$.
Then  $\mu$ is potentially isotropic. \label{cor_2019}
\end{corollary}

\emph{Proof:} Consider a sequence $\mu_1,\mu_2,\ldots$ of Borel
probability measures on $S^{n-1}$, absolutely continuous with
respect to the Lebesgue measure on $S^{n-1}$, that converges to
$\cR_* \mu$ in the weak$^*$ topology. The $\mu_i$'s are potentially
isotropic by Lemma \ref{lem_1117}. Therefore $\mu$ is potentially
isotropic according to Lemma \ref{lem_2019}. \hfill $\square$

\medskip A clever proof of Corollary \ref{cor_2019} for the
case where the measure $\mu$ is discrete and has finite support
appears in the works of Barthe \cite{barthe_inventiones}, Carlen
and Cordero-Erausquin \cite[Lemma 3.5]{CC} and Carlen, Lieb and
Loss \cite{CLL}. We were not able to
generalize their argument to the case of a general measure
satisfying (\ref{eq_1039}). The proof presented above is
unfortunately longer, but perhaps it has the advantage of being
geometrically straightforward.

\begin{lemma} Let $n$ be a positive integer and $\alpha > 0$.
Suppose that $\mu$ is an $\alpha$-decent probability measure on
$\RR^n$. Then for any $0 < \eps < 1$, there exists a linear
transformation $T: \RR^n \rightarrow \RR^n$ such that $\nu =
T_*(\mu)$ satisfies
$$
 \int_{S^{n-1}} (x \cdot \theta)^2  d \cR_*\nu(x) \leq \alpha + \eps
\ \ \ \ \ \ \ \text{for all} \ \theta \in S^{n-1}.
 $$ \label{lem_139} \end{lemma}

\emph{Proof:} By induction on the dimension $n$. The case $n=1$ is
obvious. Suppose that $n \geq 2$. We may assume that $\mu(H) < 1$
for any hyperplane $H \subset \RR^n$ that passes through the origin
(otherwise, invoke the induction hypothesis). We may also assume
that $\alpha = \sup_{E \subseteq \RR^n} \mu(E) / \dim(E)$ where the
supremum runs over all subspaces $\{ 0 \} \neq E \subseteq \RR^n$.
Corollary \ref{cor_2019} takes care of the case where
$$ \mu(E) < \dim(E) / n $$
for any subspace $E \subset \RR^n$ with $1 \leq \dim(E) \leq n-1$.
We may thus focus on the case where there exists a proper subspace
$E \subset \RR^n$ with $\mu(E) \geq \dim(E) / n$. Clearly $\alpha
\geq 1/n$. Consequently,  there is a subspace $E \subset \RR^n$,
with $1 \leq \dim(E) \leq n-1$, such that $\alpha - \eps / (3 n)
\leq \mu(E) < 1$. Let $T: \RR^n \rightarrow \RR^n$ be the map
defined by
$$ T (x) = \left \{ \begin{array}{cl} Proj_{E^{\perp}} x & x \not \in E \\ x & x \in E
\end{array} \right. $$  The map $T$
may be viewed as a ``stratified linear map'' as in Furstenberg
\cite{F}. Set $\lambda = \mu(E) > 0$. The probability measure $T_*
\mu$ on $\RR^n$ is supported on $E \cup E^{\perp}$, and it may be
decomposed as
$$ T_* \mu = \lambda \mu_E + (1 - \lambda) \mu_{E^{\perp}} $$
where $\mu_E = \mu|_E$ is the conditioning of $\mu$ to $E$, and
$\mu_{E^{\perp}}$ is a certain probability measure supported on
$E^{\perp}$. Clearly, $\mu_E = \mu|_E$ is $(\alpha /
\lambda)$-decent. Regarding $\mu_{E^{\perp}}$, let us select a
subspace $F \subseteq E^{\perp}$. Then, $\mu_{E^{\perp}}(F) = 0$ if
$F = \{ 0 \}$ and otherwise
\begin{eqnarray*}
 \lefteqn{ (1 - \lambda) \mu_{E^{\perp}}(F) = \mu(T^{-1}(F \setminus \{ 0 \})) = \mu \left( (F \oplus E) \setminus E \right)
= \mu( F \oplus E) - \mu(E) } \\ & \leq & \alpha (\dim(E) + \dim(F)) - (\alpha - \eps/(3n)) \dim(E)
\leq (\alpha + \eps/3) \dim(F)
\end{eqnarray*}
where $E \oplus F = \{ x + y ; x \in E, y \in F \}$ is the subspace
spanned by $E$ and $F$. Consequently, $\mu_{E^{\perp}}$ is an
$\left( (\alpha + \eps/3)/ (1-\lambda) \right)$-decent measure on
$E^{\perp}$. We may apply the induction hypothesis for $\mu_E$ and
$\mu_{E^{\perp}}$. We conclude that there exists a linear
transformation $S: \RR^n \rightarrow \RR^n$, with $S(E) \subseteq E$
and $S(E^{\perp}) \subseteq E^{\perp}$, such that
\begin{equation}
 \int_{S^{n-1}} (x \cdot \theta)^2 d (\cR \circ S \circ T)_* \mu(x) \leq \alpha + \eps/2
\ \ \ \ \ \text{for any} \ \theta \in S^{n-1}. \label{eq_307}
\end{equation}
The problem is that $S \circ T$ is not a linear map.
However, it is easy to approximate it by a linear map:
For $0 < \delta < 1$ denote $T_\delta x = x - \delta Proj_E x$.
Then $(\cR \circ T)(x) = \lim_{\delta \rightarrow 1^-} (\cR \circ T_{\delta})(x)$
for any $0 \neq x \in \RR^{n}$. Consequently,
$$ (\cR \circ S \circ T_{\delta})_* {\mu}
= (\cR \circ S \circ \cR \circ T_{\delta})_* {\mu} \stackrel{\delta
\rightarrow 1^-}{\longrightarrow} (\cR \circ S \circ \cR \circ T)_*
\mu = (\cR \circ S \circ T)_* \mu  $$ in the weak$^*$ topology. We
conclude that the matrices $M((\cR \circ S \circ T_{\delta})_*
{\mu})$ tend to $M( (\cR \circ S \circ T)_* \mu )$ as $\delta
\rightarrow 1^-$. Hence, by (\ref{eq_307}), for some $\delta_0 < 1$,
 $$
\int_{S^{n-1}} (x \cdot \theta)^2 d (\cR \circ S \circ T_{\delta_0} )_* \mu(x) \leq \alpha + \eps
\ \ \ \ \ \text{for any} \ \theta \in S^{n-1}. $$
The map $S \circ T_{\delta_0}$ is the desired linear transformation. This completes
the proof. \hfill $\square$

\section{Proof of the main theorem and some remarks}
\label{sec6}

\emph{Proof of Theorem \ref{main_theorem}:} Suppose that $\mu$ is an
$\eta$-decent probability measure on $\RR^n$. According to Lemma
\ref{lem_139}, there exists a linear map $S: \RR^n \rightarrow
\RR^n$ such that $\nu = S_* \mu$ satisfies
$$ \int_{S^{n-1}} (x \cdot \theta)^2 d \cR_* \nu(x) \leq 2 \eta  \ \ \ \ \ \ \text{for all} \ \ \theta \in S^{n-1}. $$
We  invoke Corollary \ref{cor_216} for the measure $\nu$. We see that
if the positive integer $d$ and $0 < \eps < 1/2$ are such that
$$ 2 \eta \leq (\tilde{c} \eps)^{\tilde{C} d}, $$
then there exists a $d \times n$ matrix $\Gamma$ for which the measure $\Gamma_* \nu$ on $\RR^d$ is $\eps$-radial proper. Setting $T = \Gamma S$, a $d \times n$ matrix, we conclude  that $T_* \mu = \Gamma_* \nu$ is a measure on $\RR^d$ which is $\eps$-radial proper. \hfill $\square$

\medskip
\emph{Proof of Corollary \ref{cor_1222}:} We may assume that $n$
exceeds a given large universal constant. Denote $d = \bar{c} \lceil
\sqrt{\log n} \rceil$ and $\delta = e^{-d} / 500$, for a small
universal constant $0 < \bar{c} < 1$ such that $n \geq (C /
\delta)^{C d}$ where $C$ is the universal constant from Theorem
\ref{main_theorem}. According to Theorem \ref{main_theorem}, we may
pass to a $d$-dimensional marginal and assume that our measure $\mu$
is a proper $\delta$-radial measure on $\RR^d$. For $t \in \RR$ and
$L > 0$ let $\chi_{t, L}$ be the $L$-Lipschitz function on the real
line which equals zero on $(-\infty, t]$ and one on $[t + 1/L,
\infty)$. Recall the Kantorovich-Rubinstein duality as in
(\ref{eq_1455}) above. Then, for any probability measure $\nu$ on
the unit sphere $S^{d-1}$ and $0 < t < c \leq 1/2$,
$$ \nu( \{ x  ; x_1 \geq t \} )
\geq \int_{S^{d-1}} \chi_{t, d}(x_1) d \nu(x) \geq
\int_{S^{d-1}} \chi_{t, d}(x_1) d \sigma_{d-1}(x) - d W_1(\nu, \sigma_{d-1}),
$$ where $x = (x_1,\ldots,x_d)$ are the coordinates of $x \in S^{d-1}$.
The integral with respect to $\sigma_{d-1}$ may be estimated
directly, and it is bounded from below by $c e^{-C t^2 d}$  (note
that the marginal of $\sigma_{d-1}$ on the first coordinate has a
density that is proportional to $(1 - t^2)_+^{(d-3)/2}$ on $[-1,1]$). We
conclude that for any $0 < t < c$ and an interval $J = [a, b]
\subset (0, \infty)$ with $\mu(S(J)) \geq \delta$,
\begin{eqnarray} \nonumber
\lefteqn{ \mu|_{S(J)} ( \{ x  ; x_1 \geq a t \} )
\geq \cR_*(\mu|_{S(J)}) ( \{ x  ; x_1 \geq t \} )} \\
& \geq & c e^{-C t^2 d} - d \cdot W_1( \cR_*(\mu|_{S(J)}), \sigma_{d-1})
\geq c e^{-C t^2 d} - d
\delta \geq c^{\prime} e^{-C^{\prime} t^2 d}.
\label{eq_1228}
\end{eqnarray}
Similarly, for any interval $J = [a, b] \subset (0, \infty)$
with $\mu(S(J)) \geq \delta$,
\begin{eqnarray} \nonumber
\lefteqn{ \mu|_{S(J)}  \left( \left \{ x ; |x_1| \geq 20 b /
\sqrt{d} \right \} \right) \leq \cR_*(\mu|_{S(J)}) \left( \left \{ x
; |x_1| \geq 20 / \sqrt{d} \right \} \right)
 } \\
& \leq & \int_{S^{d-1}} \chi_{20 /\sqrt{d} - d^{-1}, d}(|x_1|)
d \sigma_{d-1}(x) + d \cdot W_1( \cR_*(\mu|_{S(J)}), \sigma_{d-1})
 \leq 1/5, \phantom{aaaaaaa}
\label{eq_1200}
\end{eqnarray}
where the integral with respect to $\sigma_{d-1}$ is estimated in a
straightforward manner. Let $\tilde{M} >0 $ be a quantile with
$$ \mu ( \{ x ; |x| \leq \tilde{M} \}) \geq 3/4 \ \ \ \ \text{and} \
\ \ \  \mu ( \{ x ; |x| \geq \tilde{M} \}) \geq 1/4. $$ Let $a > 0$
be such that the interval $J = [a, \tilde{M}]$ satisfies $\mu(S(J))
\geq 2/3$. We apply
 (\ref{eq_1200}) for the interval $J = [a, \tilde{M}]$ to deduce that
 \begin{equation}  \mu \left( \left \{ x; |x_1| \geq 20 \tilde{M} / \sqrt{d} \right \} \right)
 \leq \frac{1}{3} +
 \frac{2}{3} \cdot  \mu|_{S(J)}  \left( \left \{ x ; |x_1| \geq 20 \tilde{M} / \sqrt{d}
\right \} \right) \leq \frac{1}{3} + \frac{2}{3} \cdot \frac{1}{5} <
\frac{1}{2}. \label{eq_1117} \end{equation}
 Suppose that $M > 0$
satisfies (\ref{eq_1159}) with the linear functional $\vphi(x)
=x_1$. We learn from (\ref{eq_1117}) that necessarily $ M \leq 20
\tilde{M} / \sqrt{d}$. Let $b > 0$ be such that the interval $J = [
\tilde{M}, b]$ satisfies $\mu(S(J)) \geq 1/5$. We apply
 (\ref{eq_1228}) for the interval $J = [\tilde{M}, b]$ and conclude
 that, for any $0 \leq t \leq c \sqrt{d} / 20$,
$$ \mu \left( \left \{ x ; x_1 \geq t M \right \} \right ) \geq \frac{1}{5} \cdot
  \mu|_{S(J)} \left(\left  \{ x ; x_1 \geq 20 t \tilde{M} / \sqrt{d} \right \}   \right)
  \geq \frac{c^{\prime}}{5}  \exp \left(
  -400 C^{\prime} t^2   \right). $$
Since $c \sqrt{d}/20 \geq \tilde{c} \log^{1/4} n$, the proof of the
lower bound for $\mu ( \{ x  ; x_1 \geq  t M \} )$ is complete. The
proof of the lower bound for $\mu ( \{ x  ; x_1 \leq -t M \} )$ is
almost entirely identical. The corollary is thus proven. \hfill
$\square$

\medskip
{\it Remarks.}

\begin{enumerate}

\item It is conceivable that a more delicate analysis yields  a better bound for $R_n$ in Corollary \ref{cor_1222}. However, note
that $R_n \leq C \sqrt{\log n}$ as is shown by the example where
$\mu$ is distributed uniformly on $n$ linearly independent vectors
in $\RR^n$. Compare the ``super-gaussian'' tail behavior of
Corollary \ref{cor_1222} with the almost sub-gaussian bounds in the
convex case in \cite{psi_two} and in Giannopoulos, Paouris and Pajor
\cite{gpp}.

\item The central limit theorem for convex bodies \cite{clt, power_law}
states that any uniform probability measure on a high-dimensional
convex set has some low-dimensional marginals that are approximately
Gaussian. It is clear that there are perfectly regular
probability measures in high dimension (e.g., a mixture of two
Gaussians) without any approximately Gaussian marginals. Therefore,
a geometric condition such as convexity is indeed relevant when we
look for approximately Gaussian marginals.
 For arbitrary high-dimensional measures without convexity properties, we may still state the more modest
conclusion that  some of the marginals are approximately
spherically-symmetric, according to Theorem \ref{main_theorem}.
There is no hope for approximate Gaussians.

\medskip Theorem \ref{main_theorem} bears
 a strong relation to the proof of  the central limit theorem for convex bodies presented in \cite{clt, power_law}
(see \cite{ptrf} for another proof, which at present works only for
a subclass of convex bodies). That proof begins by showing that
marginals of the uniform measure on a convex body are approximately
spherically-symmetric. The approximation in \cite{clt, power_law} is
rather strong compared to Theorem \ref{main_theorem}, but
nevertheless, a simple compactness argument enables us to leverage
Theorem \ref{main_theorem} in order to obtain the desired type of
approximation. In principle, this approach yields a slightly
different proof of the central limit theorem for convex sets, albeit
with weaker estimates.

\medskip The Euclidean structure with respect to
which a random projection ``works'' with high probability seems {\it
a priori} different in Theorem \ref{main_theorem} and in the central
limit theorem for convex bodies. In Theorem \ref{main_theorem} we
use the Euclidean structure with respect to which the covariance
matrix of $\cR_* \mu$ is scalar, while in the central limit theorem
for convex bodies, the most natural position is to require the
covariance matrix of $\mu$ itself to be a scalar matrix (compare
also with \cite{K_M}, \cite{em}). For convex bodies, these Euclidean
structures are close to each other, since most of the mass of a
normalized convex body is located very close to  a sphere (see
\cite{euro}).

\item The linear map $T$ in Theorem \ref{main_theorem} may be assumed
to be an orthogonal projection. This follows from the following
simple observation we learned from G. Schechtman: Any
$n$-dimensional ellipsoid has an $\lceil n/2 \rceil$-dimensional
projection which is precisely a Euclidean ball. Therefore, in order
to show that $T$ may be chosen to be an orthogonal projection, one
essentially has to verify that a $\lceil d/2 \rceil$-dimensional
marginal of an $\eps$-radial measure on $\RR^d$ is $100
\eps^{1/8}$-radial. We omit the details.

\item The isoperimetric inequality on the high-dimensional
sphere, which is the cornerstone of the concentration of measure
phenomenon (see Milman and Schechtman \cite{MS}), is not used in the
proof of Theorem \ref{main_theorem}. We do apply Levy's lemma, which
embeds the isoperimetric inequality, in the proof of Lemma
\ref{lem_456}, but only in $d$ dimensions. The dimension $d$ here is
typically not very large.

\item For a positive integer $d$ and $\eps > 0$ denote by $N_0(\eps, d)$
the minimal dimension with the following property: Whenever $N \geq
N_0(\eps, d)$, any $N$-dimensional Banach space has a
$d$-dimensional subspace which is $\eps$-close to a Hilbert space.
The classical Dvoretzky's theorem states that $N_0(\eps, d) \leq
\exp(C d / \eps^2)$, where $C
> 0$ is a universal constant (see Milman \cite{mil} and references therein).
The power of $1/\eps$ in the exponent in the bound for $N_0(\eps,d)$
can be made arbitrarily close to one at the expense of increasing
the universal constant $C$ (see Schechtman \cite{schechtman}). It is
conceivable, however, that these bounds are still far from optimal;
perhaps $N_0(\eps,d)$ can be made as small as $(C / \eps)^{C d}$?
See Milman \cite{few} for a discussion of this conjecture. An
affirmative answer for the case $d=2$ was given by Gromov
\cite{few}, using a topological argument which does not seem to
generalize to higher dimensions.

\medskip The analogy with the present article suggests
to try and use Theorem \ref{main_theorem}, or ideas from its
proof, in order to improve the bounds in Dvoretzky's theorem.
Furthermore, the operation of marginal is dual, via the Fourier
transform, to the operation of restriction to a subspace. So, for
instance, suppose a norm $\| \cdot \|$ in $\RR^n$ may be
represented as
\begin{equation}
 \| x \| = \int_{\RR^n} |x \cdot \theta| d \mu(\theta) \label{eq_1123}
 \end{equation}
for a compactly-supported probability measure $\mu$ on $\RR^n$. In
this case, we may consider subspaces $E \subset \RR^n$ for which
$(Proj_E)_* \mu$ is $\eps$-radial, and expect that the restriction
of $\| \cdot \|$ to these subspaces is close, in a certain sense, to
the Euclidean norm. See Koldobsky \cite[Chapter 6]{kold} for a
comprehensive discussion of norms admitting representations in the
spirit of (\ref{eq_1123}).

\medskip While this approach may possibly yield some meaningful
estimates for some classes of normed spaces, it has limitations.
Theorem \ref{main_theorem} is proven by considering a {\it random}
marginal with respect to an appropriate Euclidean structure, i.e., a
projection of the given measure to a subspace which is distributed
uniformly over the Grassmannian of all $d$-dimensional subspaces in
$\RR^n$. However, for Banach spaces such as $\ell_{\infty}^N$, a
random subspace is not sufficiently close to a Hilbert space (see
Schechtman \cite{schechtman2}), and there are better choices than
the random one. (Indeed, the $\ell_{\infty}^N$ norm cannot be
represented as in (\ref{eq_1123}) or in a similar way, see Theorem
6.13 in Koldobsky \cite{kold}, going back to Misiewicz). A direct
application of Theorem \ref{main_theorem} is thus quite unlikely to
provide new information regarding approximately Hilbertian subspaces
for  {\it all} finite-dimensional normed spaces.

\item In principle, the measures $T_*(\mu)$ in
Theorem \ref{main_theorem} are not only approximately radial, but
are also approximately a composition of isotropic Gaussians. Indeed,
it is well-known that any $d$-dimensional marginal of the measure
$\sigma_{k-1}$, for $d \ll k$, is approximately an isotropic
$d$-dimensional Gaussian measure. Thus, we may project an
approximately-radial measure on $\RR^k$ to any $d$-dimensional
subspace, and obtain a measure which is approximately, in some
sense, a composition of isotropic Gaussians. We did not rigorously
investigate
this approximation property on a precise, quantitative level.

\end{enumerate}

\section{Infinite-dimensional spaces}
\label{sec7}

This section contains a corollary to Theorem \ref{main_theorem},
pertaining to probability measures supported on infinite-dimensional
spaces. We begin with a lemma regarding distributions on finite-dimensional spaces.
Let $n \geq 1$ be an integer, suppose that $\mu$
is a Borel probability measure on $\RR^n$ and let $0 < a \leq 1$. A
subspace $E \subseteq \RR^n$ is ``$a$-basic for $\mu$'' if
\begin{enumerate}
 \item[(i)] $\mu(E) \geq  a$
\item[(ii)] $\mu(F) < a$ for any proper subspace $F \subsetneq E$.
\end{enumerate}
Note that any subspace $E \subseteq \RR^n$ with $\mu(E) \geq a$ contains
an $a$-basic subspace. Also, suppose $T: \RR^n \rightarrow \RR^m$ is a linear map,
and let $E \subseteq \RR^n$ be an $a$-basic subspace for $\mu$ containing the kernel of $T$.
Then $T(E)$ is $a$-basic for $T_*(\mu)$.

\begin{lemma} Let $n \geq 1$ be an integer, $0 < a \leq 1$,
and let $\mu$ be a Borel probability measure on $\RR^n$. Then there are only finitely
many subspaces $E \subseteq \RR^n$ that are $a$-basic for $\mu$.
\label{lem_1212}
\end{lemma}

\emph{Proof:} Let $k \geq 0$ be an integer and $0 < a \leq 1$. We
will prove by induction on $k$ the following statement: For any
integer $n \geq k$ and for any Borel probability measure $\mu$ on
$\RR^n$, there are at most finitely many subspaces $E \subseteq \RR^n$
whose dimension is at most $k$, that are  $a$-basic for the measure
$\mu$. The statement clearly implies the lemma. The case $k = 0$ is
easy, as there is only one $0$-dimensional subspace in $\RR^n$.

\medskip Let $k \geq 1$. Suppose that $n \geq k$ is an integer, $0 <
a \leq 1$ and let $\mu$ be a Borel probability measure on $\RR^n$.
Denote by $\cG$ the family of all subspaces $E \subseteq \RR^n$ whose
dimension is at most $k$ that are $a$-basic for the measure $\mu$.
We need to show that
\begin{equation}
\#(\cG) < \infty. \label{eq_11200}
\end{equation}
First, note that it is sufficient to prove (\ref{eq_11200}) under the
additional assumption that $\mu ( \{ 0 \} ) = 0$. Indeed, denote
$\eps = \mu ( \{ 0 \} )$. If $a \leq \eps$ then there is only one
$a$-basic subspace in $\RR^n$, which is the subspace $\{ 0 \}$,
and (\ref{eq_11200}) clearly holds. In the non-trivial case where $a > \eps$,
we may replace $\mu$ by $\left( \mu - \eps \delta_0 \right) / (1 -
\eps)$ and $a$ by $\left( a - \eps \right) / ( 1 - \eps)$. The
family of basic subspaces remains exactly the same. From now on,
we will thus assume that $\mu ( \{ 0 \}) = 0$.

\medskip Denote by $\cE \subseteq \cG$ the collection of all subspaces $E \subseteq \RR^n$
that are $a$-basic for $\mu$, with $\dim(E) \leq k$, for which
$\mu(F) < a^2 / 8$ for any proper subspace $F \subsetneq E$. We will
prove that
\begin{equation}
 \#(\cE) \leq 2 / a < \infty.
\label{eq_11230}
\end{equation}
To that end, let $\tilde{\cE}$ be any finite subset of $\cE$, and
denote $N = \#(\tilde{\cE})$. For any two distinct subspaces $E_1,
E_2 \in \tilde{\cE}$, we have $\mu(E_1 \cap E_2) < a^2 / 8$ as $E_1
\cap E_2$ is a proper subspace of $E_1$. According to the
inclusion-exclusion principle,
$$ 1 \geq \mu \left( \bigcup_{E \in \tilde{\cE}} E \right) \geq \sum_{E \in \tilde{\cE}} \mu(E) - \sum_{E_1, E_2 \in \tilde{\cE} \atop{E_1 \neq E_2}} \mu(E_1 \cap E_2) > N a - \frac{N(N-1)}{2} a^2/8, $$
where we used the fact that $\mu(E) \geq a$ for any $E \in
\tilde{\cE}$, since $E$ is $a$-basic. We conclude that
\begin{equation}
 1 > N a - \frac{N(N-1) a^2}{16} \geq N a \left[ 1 - \frac{N a}{10} \right] \ \ \ \Longrightarrow \ \ \
 |Na -5| > 3 .
\label{eq_1136}
\end{equation}
Thus, there are no finite subsets of $\cE$ whose cardinality is $N =
\lceil 2 / a \rceil$: In this case $2 \leq N a \leq 3 $ which is
impossible according to (\ref{eq_1136}). Hence $\#(\cE) \leq 2 / a$
and (\ref{eq_11230}) is proven.

\medskip Next, denote by $\tilde{\cG}$ the family of all subspaces $E \subseteq \RR^n$
that are $a$-basic, with $\dim(E) \leq k$, for which there exists a
proper subspace $F \subsetneq E$ with $\mu(F) \geq a^2 / 8$. In view
of (\ref{eq_11230}), in order to deduce (\ref{eq_11200})  it
suffices to show that
\begin{equation}
 \#(\tilde{\cG}) < \infty.
\label{eq_1120}
\end{equation}
Whenever a subspace $E \subseteq \RR^n$ contains a proper subspace
$F \subsetneq E$ with $\mu(F) \geq a^2 / 8$, it also contains an
$a^2/8$-basic proper subspace $\tilde{F} \subsetneq E$ with
$\dim(\tilde{F}) \leq k-1$.
 By the induction hypothesis, there are only finitely many
 subspaces $\tilde{F} \subseteq \RR^n$ that are $a^2/8$-basic for $\mu$ whose dimension is at most $k-1$.
Fix such an $a^2/8$-basic subspace $\tilde{F}$. Let $\cF$ be the
collection of all subspaces $E \subseteq \RR^n$ that are
$a$-basic, contain $\tilde{F}$, and satisfy $\dim(E) \leq k$. The
task of proving (\ref{eq_1120}) and completing the proof of the
lemma is reduced to showing that
$$
 \#(\cF) < \infty. $$
Note that $\dim(\tilde{F}) \geq 1 $ as $\mu ( \{ 0 \} ) = 0  < a^2 /
8$, and hence $\{ 0 \}$ is not an $a^2/8$-basic subspace. Denote
by $P = Proj_{\tilde{F}^{\perp}}$ the orthogonal projection operator
onto $\tilde{F}^{\perp}$ in $\RR^n$. Then $\nu = P_*(\mu)$ is a
Borel probability measure on $\tilde{F}^{\perp}$. For any $E \in
\cF$, the subspace $P(E)$ is an $a$-basic subspace for the measure
$\nu$, and $\dim(P(E)) = \dim(E) - \dim(\tilde{F}) \leq k-1$. From
the induction hypothesis, we see that the set $\{ P(E) ; E \in \cF
\}$ is finite. However, $P(E_1) \neq P(E_2)$ for any distinct $E_1,
E_2 \in \cF$. Thus $\#(\cF) <\infty$, as promised. The lemma is
proven. \hfill $\square$

\medskip An alternative proof of Lemma \ref{lem_1212} was suggested by N. Alon.
His idea is to replace the first part of the proof
of the induction step with the known fact that there exists a finite set $A \subset \RR^n$
that intersects any subspace of measure at least $a$ (see, e.g., Alon and Spencer \cite[Section 13.4]{AS}).

\medskip
We write $\RR^{\infty}$ for the linear space of infinite
sequences $a = (a_1,a_2,\ldots)$
with $a_i \in \RR$ for all $i \geq 1$. The space $\RR^{\infty}$ is endowed
with the standard product topology (also known as Tychonoff's topology) and
the corresponding Borel $\sigma$-algebra.
 The projection map
$P_n: \RR^{\infty} \rightarrow \RR^n$ is defined by
$$ P_n(x) = (x_1,\ldots,x_n) $$
for $x = (x_1,x_2,\ldots) \in \RR^{\infty}$. Then $P_n$ is a continuous,
linear map. Note that any finite-dimensional subspace $E \subset \RR^{\infty}$
is a closed set. Also for any subspace $E \subseteq \RR^{\infty}$ we have
\begin{equation}
 \dim(E) = \sup_n \dim(P_n(E)).
\label{eq_1021}
\end{equation}
With a slight abuse of notation, for $m \geq n \geq 1$ we also write
$P_n : \RR^m \rightarrow \RR^n$ for the projection operator defined
by $P_n(x_1,\ldots,x_m) = (x_1,\ldots,x_n)$. We will also use the
ridiculous space $\RR^0 = \{ 0 \}$, and $P_0(x) = 0$ for any $x$.
Let $\eps > 0$ and let $X$ be a measurable linear space in which all
finite-dimensional subspaces are measurable.
A probability measure $\mu$ on $X$ is called $\eps$-decent if
for any finite-dimensional subspace $E \subseteq X$,
$$ \mu(E) \leq \eps \dim(E). $$

\begin{lemma} Let $\eps > 0$ and let $\mu$ be a Borel probability measure on
$\RR^{\infty}$. Suppose that $\mu$ is $\eps$-decent. Then, there
exists $N \geq 1$ such that $(P_N)_* \mu$ is $2 \eps$-decent.
\label{lem_2223}
\end{lemma}

\emph{Proof:} For $n \geq 0$ denote $\mu_n = (P_n)_* \mu$,
a Borel probability measure on $\RR^n$.  We say that a
subspace $E \subseteq \RR^n$ is ``thick'' if $\mu_n(E) \geq 2\eps \dim(E)$.
A thick subspace $E \subseteq \RR^n$ is necessarily of dimension at most $(2 \eps)^{-1}$.
We say that $E$ is a ``primitive, thick subspace'' if it is thick
and additionally
$$ \mu_n(F) < 2 \eps \dim(F) $$
for any proper subspace $F \subsetneq E$. Clearly, any thick subspace $E \subseteq \RR^n$
contains a primitive, thick subspace.
Observe also that a primitive, thick, $k$-dimensional
subspace $E \subseteq \RR^n$ is necessarily $2 \eps k$-basic for the
measure $\mu_n$.
From Lemma \ref{lem_1212} we thus learn that for any $n$,
 there are only finitely many primitive, thick subspaces $E \subseteq \RR^n$.

\medskip Denote by $\cV$ the collection
of all pairs $(E,n)$ such that $E \subseteq \RR^n$ is a primitive, thick subspace.
In order to prove the lemma, it suffices to show that $\cV$ is finite. Indeed,
in this case, set $N = \max \{ n + 1; \exists E \subseteq \RR^n, (E, n) \in \cV \}$.
Then there are no primitive, thick subspaces in $\RR^N$, and hence there are no
thick subspaces in $\RR^N$. Consequently $\mu_N = (P_N)_*(\mu)$ is $(2 \eps)$-decent,
and the lemma is proven. The rest of the argument is
thus concerned with the proof that $\cV$ is finite.

\medskip Define a directed graph structure on $\cV$ as follows:
There is an edge going from  the node $(E,n) \in \cV$ to the node $(F, n+1) \in \cV$
if and only if $E \subseteq P_n(F)$. Note that for each node $(F, n+1) \in \cV$,
the subspace $P_{n}(F) \subseteq \RR^{n}$ is clearly thick, hence it contains a
primitive, thick subspace $\tilde{E} \subseteq \RR^n$. Therefore each node $(F, n+1)$
is connected to a certain node $(\tilde{E}, n) \in \cV$. We conclude that there is a
path from $(\{ 0 \}, 0) \in \cV$ to any node in $\cV$.  For each $n \geq 1$ there
are only finitely many nodes of the form $(E, n) \in \cV$, since
 there are only finitely many primitive, thick subspaces $E \subseteq \RR^n$.
 Therefore, $\cV$ is finite
if and only if it does not contain an infinite path.

\medskip We deduce that in order to prove the lemma, it suffices to show that there is no
sequences of subspaces $E_n \subseteq \RR^n \ (n=0,1,\ldots)$
such that for any $n \geq 0$,
\begin{equation}
 E_n \subseteq P_n(E_{n+1}) \ \ \ \text{and} \ \ \ (E_n, n) \in \cV.
\label{eq_1009}
\end{equation}
Assume by contradiction that a sequence of subspaces satisfying (\ref{eq_1009})
exists. Recall that a subspace of dimension larger than $(2 \eps)^{-1}$ cannot be thick,
hence $\dim(E_n)$ is bounded by $(2 \eps)^{-1}$. Additionally,
 $\dim(E_n) \leq \dim(E_{n+1})$ for all $n$.
Therefore, there exist $n_0 \geq 1$ and $d \leq (2 \eps)^{-1}$ such that
$$ \dim(E_n) = d \ \ \ \ \ \text{for all} \ n \geq n_0. $$
Consequently, $E_n = P_n(E_{n+1})$ for any $n \geq n_0$.
Consider the direct limit
$$ E = \left \{ a \in \RR^{\infty} \, ; \, P_n(a) \in E_n \ \text{for all} \ n \geq n_0 \right \} \subseteq \RR^{\infty}. $$
Then $E = \cap_{n \geq n_0} P_n^{-1}(E_n)$ is a subspace of $\RR^{\infty}$
with $P_n(E) = E_n$ for all $n \geq n_0$.
Furthermore, $\dim(E) = d$ according to (\ref{eq_1021}). Note that
$P_n^{-1}(E_n) \supset P_{n+1}^{-1}(E_{n+1})$ for any $n \geq n_0$. Therefore
$$ \mu(E) = \mu \left( \bigcap_{n \geq n_0} P_n^{-1}(E_n) \right)
= \lim_{n \rightarrow \infty} \mu \left( P_n^{-1}(E_n) \right) =
\lim_{n \rightarrow \infty} \mu_n \left( E_n \right) \geq 2 \eps d,
$$ since $E_n \subseteq \RR^n$ is a $d$-dimensional thick subspace.
Hence $\mu(E) \geq 2 \eps d$, in contradiction to our assumption
that $\mu$ is $\eps$-decent. We conclude that there are no infinite paths in
$\cV$, and hence that $\cV$ is finite .The lemma is proven. \hfill
$\square$

\medskip Suppose $X$ is a topological vector space. We say that $X$
has a countable separating family of continuous, linear functionals
if there exist continuous linear functionals $f_1,f_2,\ldots: X
\rightarrow \RR$ such that for any $x \in X$,
$$ x = 0 \ \ \ \ \ \Longleftrightarrow \ \ \ \ \ \ \forall n, \
f_n(x) = 0. $$ This condition is not too restrictive. For example,
any separable normed space, any separable Fr\'echet space, and any
 topological vector space dual to a separable Fr\'echet space --
admits a countable, separating family of continuous, linear
functionals.

\begin{corollary} Let $\eps > 0$, let $d \geq 1$ be an integer,
and let $X$ be a topological vector space with a countable separating
family of continuous, linear functionals.
Suppose that $\mu$ is an $\eps$-decent Borel probability measure on
$X$. Then, there exists a continuous linear map $T: X \rightarrow
\RR^d$ such that $T_*(\mu)$ is $\delta$-radial proper, for $\delta =
c \eps^{c/d}$. Here, $c > 0$ is a universal constant.
\label{cor_1251}
\end{corollary}

\emph{Proof:} Let $f_1,f_2,\ldots: X \rightarrow \RR$ be the
separating sequence of continuous linear functionals. Then the
linear map $T: X \rightarrow \RR^{\infty}$ defined by
$$ T(x) = (f_1(x),f_2(x),\ldots)  $$
is a continuous linear embedding. Since $\mu$ is $\eps$-decent, then
also $T_*(\mu)$ is an $\eps$-decent, Borel probability measure on
$\RR^{\infty}$. According to Lemma \ref{lem_2223}, there exists a
finite $N \geq 1$ and a continous linear map $P: \RR^{\infty} \rightarrow
\RR^N$ such that $(P \circ T)_*(\mu)$ is a $2 \eps$-decent measure
on $\RR^N$. The corollary now follows from Theorem
\ref{main_theorem}. \hfill $\square$

\medskip Note that
the linear map $T$ in Corollary \ref{cor_1251} is not only measurable
but also continuous. In
principle, we could have formulated Corollary \ref{cor_1251} for a
probability measure on a measurable linear space, without having to
rely on an ambient topology: All we need is a linear, measurable
embedding in $\RR^{\infty}$. We refer the reader to Tsirelson
\cite{Ts} for a discussion of measures on infinite-dimensional linear spaces, and
for an exposition of Vershik's ``de-topologization'' program \cite{V1, V2}.
We conclude this note with an infinite-dimensional
analog of Corollary \ref{cor_1222}.

\begin{corollary} Let $X$ be a topological vector space with a countable separating
family of continuous, linear functionals.
Suppose that $\mu$ is a Borel probability measure on
$X$ such that $ \mu(E) = 0$ for any finite-dimensional subspace $E \subset X$.
Then, for any $R > 0$, there exists a non-zero, continuous linear functional
$\vphi: X \rightarrow \RR$ such that
$$ \mu \left( \{ x  ; \vphi(x) \geq t M \right \}) \geq c \exp(-C t^2)
\ \ \ \ \ \ \text{for all} \ 0 \leq t \leq R $$ and
$$ \mu \left( \{ x  ; \vphi(x) \leq -t M \right \}) \geq c \exp(-C t^2)
\ \ \ \ \ \ \text{for all} \ 0 \leq t \leq R $$ where $M > 0$ is a
median, that is,
$$
 \mu \left( \{ x  ; |\vphi(x)| \leq M \right \})
\geq 1/2 \ \ \ \ \text{and} \ \ \ \ \mu \left( \{ x  ; |\vphi(x)|
\geq M \right \}) \geq 1/2 $$ and $c, C
> 0$ are universal constants.
\end{corollary}

\bigskip

{\small \noindent  School of Mathematical Sciences, Tel-Aviv
University, Tel-Aviv 69978, Israel

{\small \noindent \it e-mail address:} {\small \verb"klartagb@post.tau.ac.il"}

\end{document}